 \documentclass[12pt]{amsart}
\usepackage{amsmath}
\usepackage{enumerate}
\usepackage{amssymb}
 
\usepackage{amsbsy}
\usepackage{amsfonts}
 \usepackage{tikz}
\newtheorem{theorem}{Theorem}

\newtheorem{lemma}{Lemma}

\def\min{\mbox{min}}
\newcommand{\sm}{\left(\begin{smallmatrix}}
\newcommand{\esm}{\end{smallmatrix}\right)}
\theoremstyle{definition}

  \DeclareMathAlphabet{\newcal}{U}{dutchcal}{m}{n}
  \DeclareMathAlphabet{\eucal}{U}{eus}{m}{n}

 \def\R{\mathbb R} \def\C{\mathbb C} \def\Z{\mathbb Z}  
\def\={\;=\;} \def\+{\,+\,} \def\m{\,-\,}  \def\df{\;:=\;}  \def\la{\langle} \def\ra{\rangle}  
\def\e{\varepsilon}  \def\t{\tau}  \def\g{\gamma} \def\G{\Gamma}     \def\j{V}
      
\def\CC{\mathbf C}      \def\Cg{\mathfrak C}    \def\Rg{\mathfrak R}   
\def\BB{\mathcal B}      
    
\def\be{\begin{equation}}  \def\ee{\end{equation}}  \def\bes{\begin{equation*}}  \def\ees{\end{equation*}}
\def\ba{\begin{align}}  \def\bas{\begin{align*}}  \def\ea{\end{align}}  \def\eas{\end{align*}}
\theoremstyle{remark}

 \usepackage{tikz}

\DeclareMathAlphabet{\newcal}{U}{dutchcal}{m}{n}
  \DeclareMathAlphabet{\eucal}{U}{eus}{m}{n}

 \def\R{\mathbb R} \def\C{\mathbb C} \def\Z{\mathbb Z}  
\def\={\;=\;} \def\+{\,+\,} \def\m{\,-\,}  \def\df{\;:=\;}  \def\la{\langle} \def\ra{\rangle}  
\def\e{\varepsilon}  \def\t{\tau}  \def\g{\gamma} \def\G{\Gamma}     \def\j{V}
      
\def\CC{\mathbf C}      \def\Cg{\mathfrak C}    \def\Rg{\mathfrak R}   
\def\BB{\mathcal B}      
    
\def\be{\begin{equation}}  \def\ee{\end{equation}}  \def\bes{\begin{equation*}}  \def\ees{\end{equation*}}
\def\ba{\begin{align}}  \def\bas{\begin{align*}}  \def\ea{\end{align}}  \def\eas{\end{align*}}

\theoremstyle{remark}

\theoremstyle{plain}
\newtheorem{thm}{Theorem}[section]
\newtheorem{prop}[thm]{Proposition}

\theoremstyle{definition}

\newtheorem{rmk}[thm]{Remark}
\numberwithin{equation}{section}

\newfont{\FieldFont}{msbm10 scaled\magstep1}

\newcommand{\pf}{\noindent\bf Proof }

\def\C{\mathbb C}

\numberwithin{equation}{section}

\newcommand{\GL}{\mbox{GL}}

\def\level{N}
\def\divlevel{M}

\usepackage{color}

\definecolor{blue}{rgb}{0,0,1}
\definecolor{red}{rgb}{1,0,0}
\definecolor{green}{rgb}{0,.6,.2}
\definecolor{purple}{rgb}{1,0,1}
\long\def\red#1\endred{{\color{red}#1}}

\long\def\blue#1\endblue{{\color{blue}#1}}
\long\def\purple#1\endpurple{{\color{purple}#1}}
\long\def\green#1\endgreen{{\color{green}#1}}

\begin{document}

\title{Twisted Kronecker Series and Periods of Modular forms on $\Gamma_0(N)$   }

 \author{Clifford Blakestad }
\author{YoungJu Choie  }
\address{Clifford Blakestad \endgraf
}
\email{math.blakestad@gmail.com  }

\address{YoungJu Choie \endgraf
Department of Mathematics\endgraf
Pohang University of Science and Technology \endgraf
Pohang,
Republic of Korea
}
\email{yjc@postech.ac.kr}

   \thanks{This work was partially supported by  NRF-2022R1A2B5B0100187113 and   2021R1A6A1A10042944  }

 
\maketitle
\maketitle

\begin{abstract} 
We introduce  an infinite family of Kronecker series twisted by  characters.
As an application,   we give a closed formula for the sum of  all Hecke eigenforms on $\Gamma_0(N)$ multiplied  by their twisted period polynomials in terms of the product of those twisted Kronecker series, when $N$ is square free.
This extends  an identity of Zagier  among period polynomials, Hecke eigenforms and  a  quotient of Jacobi theta series. \\
 \end{abstract}

\bigskip

\section{{\bf{Introduction and Statement of Main Theorem }}}
 Let us recall the following Kronecker series, for $\t \in \mathbb{H}, u, v\in \mathbb{C},$  
\begin{equation}\label{kro}
F_{\t}(u,v)=\sum_{n\in \mathbb{Z}} \frac{\eta^n}{q^n\xi-1},  \, \,   q=e^{2\pi i \t}, \xi=e^u, \,\eta=e^v,\,  
  |q|<|\xi|<1   
\end{equation}
 introduced by Kronecker and studied by Weil \cite{Weil} who showed that $F_{\t}(u,v)$ can be expressed as
 $$F_{\t}(u,v)=\frac{\theta'(0) \theta(u+v)}{\theta(u) \theta(v)}$$  with Jacobi theta series
$$\theta(u)=q^{\frac{1}{8}}(\xi^{\frac{1}{2}}-\xi^{-\frac{1}{2}})
\prod_{n\geq 1} (1-q^n)(1-q^n\xi)(1-q^n\xi^{-1}).$$
It turns out that this Kronecker series appears in   a vast array of research areas.
For instance, it shows up in topology when studying
local elliptic classes of Bott-Samelson varieties \cite{RRW},  
in  combinatorics when counting  torus covers, when studying 
    characters of topological vertex algebra \cite{HE}  and  
  in various other  places such as in   mathematical physics \cite{AW, CMZ} and  in the theory of  modular forms \cite{ BFI, F, DR,  pp, JMO}.
Of particular interest for this paper,  Zagier \cite{Z1991} found a striking identity between a generating  function   whose coefficients encode Hecke eigenforms together with their periods and  the above Kronecker series. \\

Let us recall the period polynomial of an elliptic  cusp   form  
$f(\tau)=\sum_{\ell\geq 1 }a_f(\ell)q^{\ell}\, \,  (\tau \in \mathbb{H}=$  upper half plane, $q=e^{2\pi i \tau})$ of weight $k$ on $SL_2(\mathbb{Z})$ is the polynomial of degree $k-2$ defined by
 \begin{eqnarray*}\label{pe}r_f(X)=\int_{0}^{i\infty}f(\tau)(\tau-X)^{k-2}d\tau,\end{eqnarray*}
 \noindent or equivalently  by
$$r_f(X)=-\sum_{n=0}^{k-2}\frac{(k-2)!}{(k-2-n)!} \frac{L(f,n+1)}{(2\pi i)^{n+1}} X^{k-2-n} 
$$
 where  $L(f,s)=\sum_{n=1}^{\infty}\frac{a_f(n)}{n^s} \,  (re(s)\gg 0) .$
The maps $f \mapsto r_f^{\rm ev}$ and $f \mapsto  r_f^{\rm od}$ assigning to $f$ the even and odd parts of $r_f$ are both injective into the ring of two-variable polynomials,  with known images  from the  Eichler-Shimura-Manin  and Kohnen-Zagier theory \cite{KZ1982}. 
 Define the two-variable polynomials 
\begin{eqnarray*}
\Rg_f(X,Y)&:=&\frac{r_f^{\rm ev}(X)r_f^{\rm od}(Y)  }{(2i)^{k-3} <f, \, f>} \\
R_f(X,Y)&:=&\Rg_f(X,Y)+\Rg_f(Y,X).
\end{eqnarray*}
 In the case that $f$ is a Hecke eigenform it turns out that $\Rg_f(X,Y)\in \mathbb{Q }_f [X,Y]$,   where $ \mathbb{Q}_f$ is the field generated by Fourier coefficients of $f$ over $\mathbb{Q}.$
\medskip

  Zagier's   identity relates the following   generating function containing  all Hecke eigenforms  together  with all   critical values, to the Kronecker series $F_{\tau}(u,v)$ studied in Weil \cite{Weil} : 

\begin{eqnarray}\label{Za1991}
 &&   \frac{(XY-1)(X+Y)}{X^2Y^2}T^{-2}
  +\sum_{k\geq 2}
\sum_{  f \in \mathcal{B}_k
 }
R_f(X,Y) f(\tau) \frac{T^{k-2}}{(k-2)!}
\\
&=&
F_{\tau}(XT, YT)  F_{\tau}( T, -XYT) ,  \nonumber
 \end{eqnarray}

\noindent where $\mathcal{B}_k$ is a set of all Hecke eigenforms of weight $k$ on $SL_2(\mathbb{Z})$.
 The above  expansion  with respect to the variable  $ T $   gives an algorithm to compute Hecke  eigenforms (see more details in \cite{Z1991}).  
  It turns out that   identities   such as   (\ref{Za1991})  are   not accidental but
   in fact exist  for a general group  $\Gamma_0(N)$ (see \cite{CPZ}) as well as   for  Hilbert modular forms \cite{C-hilbert}.  It took a long time to extend the existence of such an identity to  more general   groups due to  the  complexity of the structure  of the associated spaces of modular forms.  
\medskip

  At the same time,  the twisted  $L$-functions
 $L(f,\chi,s)=\sum_{n\geq 1}\frac{\chi(n) a_f(n)}{n^s}$
associated   to   cusp forms   $f$   twisted   by a character $\chi$  
play  an important role in number theory. 
For instance,   Waldspurger showed that  the algebraic part of central values of $L(f,\chi,s)$ are    squares (see  \cite{KZ1982}).
 The   Goldfeld conjecture \cite{Go} predicts the behaviour of the   average sum of the order of zeros of 
$L(f, \chi, s),$  which   according to the Birch and Swinnerton-Dyer conjecture  is the  average sum of the rank of the group of rational points $E_{\chi}(\mathbb{Q})$ on twisted elliptic curve s  $E_{\chi},$ when $f$ is an associated weight $2$ newform to $E.$
 \medskip

In this paper 
we  introduce an infinite family of Kronecker series   twisted by characters $F^{\chi}_{\tau}(u,v)$   and 
extend Zagier's identity   (\ref{Za1991}) to  one    between a
generating function for twisted periods of    modular forms    and  those  twisted     Kronecker  series (see Theorem \ref{main}).   
Despite a long history and a vast application  of  Kronecker series   in term of many different   perspectives (see \cite{BK, BC, C-hilbert}),   it is   surprising   that twisted   Kronecker series  have not been studied before.  

This work is  an   attempt to have an arithmetic application of the   introduced   twisted Kronecker series.  
  The authors also  further explore   some of the  geometric  and  arithmetic properties of the twisted Kronecker series   in   upcoming work \cite{BC}. 
\medskip

 As in the classical case, 
for a Hecke eigenform  $f(\tau)=\sum_{n\geq 1} a_f(n)q^{  n },  q=e^{2\pi i \tau},   $  of weight $k$ on $\G_0(N),$ 
the {\it period polynomial} $r_f(X)$ of $f$ is  a
polynomial  of degree $k-2$ defined by
 
 \begin{eqnarray*}
r_{f }(X)
 &=&  \sum_{n=0}^{k-2}  (-1)^n  \frac{(k-2)!}{  n!   (k-2-n)!  }   
r_n(f ) X^{k-2-n} , 
\end{eqnarray*}
where 
\begin{eqnarray*}
r_{n}(f ) =  \frac{  i^{n+1}   \Gamma(n+1)}{(2\pi  )^{n+1} }  L(f, n+1), 
L(f,s):=\sum_{n\geq 1}\frac{a_f(n)}{n^s} \, \,  ( \, re(s)\gg 0).
\end{eqnarray*}
Further if  $f|_kW_N=\epsilon(N) f $   then the functional equation of $L(f,s)$ implies that 
\begin{eqnarray}\label{fun1}
r_{k-n-2}(f)= (-1)^{n+1}\epsilon(N) 
N^{ - \frac{k}{2}+n+1} r_n(f).
\end{eqnarray} 
For twists   
$f_{\chi}(\tau)=\sum_{n\geq 1}^{\infty}\chi(n) a_f(n)q^{  n }\in  S_{k}(\Gamma_0(N^2), \chi^2) $ 
by a Dirichlet character with conductor $N$,
the   {\it{twisted period polynomial}} $r_{f_{\chi}}(X)$  is similarly  given by 
  \begin{eqnarray*}
r_{f_{\chi}}(X)
 &=&  \sum_{n=0}^{k-2}  (-1)^{n } \frac{(k-2)!}{  n!  (k-2-n)!  }   
r_n(f_{\chi}) X^{k-2-n},
\end{eqnarray*}
where 
\begin{eqnarray}\label{tt}
&&  \\
&& r_{n}(f_{\chi} ) =  \frac{  i^{n+1}\Gamma(n+1)}{(2\pi )^{n+1} }  L(f,\chi,n+1), \, L(f, \chi.s)=\sum_{n\geq 1} \frac{\chi(n) a_f(n)}{n^s}\,\,( \, re(s) \gg 0).
 \nonumber
 \end{eqnarray}
  Since  $ f_{\chi}|_kW_{N^2}=\lambda_N f_{\overline{\chi}} $ with
$\lambda_N=\chi(-1)\frac{W(\chi)}{W(\overline{\chi})}$   (For $\lambda_p, p|N,$ see Theorem 4.1 in \cite{AL} .  See Proposition \ref{tfun} in Section 3.1 for details)
so  the functional equation of $L(f, \chi, s)$  gives     

 
\begin{eqnarray}\label{twist-fun}
 r_{k-2-n}(f_\chi)=(-1)^{n+1} \chi(-1)\frac{W(\chi)}{W(\overline{\chi})}N^{2n+2-k}r_{n}(f_{\overline{\chi}}),
\end{eqnarray}  
 where $W(\chi)$ is the Gauss sum (see 2.1 Notations).

\noindent Write $R_{f_{\chi}} (X,Y)$ as $\Rg_{f_{\chi}}(X,Y)+\Rg_{ {f_{\chi}}}(Y,X)$ with

\begin{eqnarray}\label{pp}
&& 
\Rg_{f_{\chi}}(X,Y) 
=\frac{1}{2}\bigl(\Cg_{f_{\chi}}(X,Y)+(XY)^{k-2}\Cg_{f_{\chi}}(-\frac{1}{X}, -\frac{1}{Y})\bigr) \nonumber \\
&& \Cg_{f_{\chi}}(X,Y)
= \frac{ 
   r^{\rm ev}_f(\frac{ Y}{N}) \, \, r^{\rm od}_{f_{\chi}}( \frac{X}{N } ) +    r^{\rm ev}_{f_{\chi}}(\frac{Y}{N})\,\, r^{\rm od}_{f }(\frac{X}{N})   }{N^{1-k}\,\, W(\chi) 2(2 i)^{k-3} \,\langle f,\, f\rangle}.
\end{eqnarray}
This two variable polynomial $R_{f_{\chi}}(X,Y)$ transforms under $\sigma\in Gal(\mathbb{C} / \mathbb{Q})$
by $R_{\sigma(f _{\chi})}=\sigma(R_{f_{\chi}}), $ so $R_{f_{\chi}}$
 has  coefficients in the number field $\mathbb{Q}_f(\chi)$
 (see Theorem \ref{Razar} in  Section 3 ).   
\bigskip

\begin{theorem}\label{main} For every integer~$k\ge2$ and an even primitive character $\chi$ of square-free conductor $N \geq 1 $, define
$$\label{Main1}  
{\bf{C}}_{k,N, \chi}(X,Y,\t) \df \frac1{(k-2)!}\,\sum_{f\in\BB_{k,N}} R_{f_{\chi} }(X,Y) \,f(\t)\,,$$
where $\mathcal{B}_{k,N}$   defined in Section 2  
is the basis of Hecke forms on~$\G_0(N)$.  Then the generating function
\begin{eqnarray*}\label{Defgenfn} 
\CC_{N,\chi}(X,Y,\t,T)  \df  & &   
  \chi(0) \frac{ (X+Y)( XY-1)}{ X^2Y^2T^2}   \\
&+& \sum_{k\geq 2}\CC_{k, N, \chi}(X,Y,\t)\,   T^{k-2} \in \mathbb{Q}(\chi)[X,Y][[q,T]]
\end{eqnarray*} 
is given in terms of   the twisted Kronecker series $F_{\tau}^{\chi}( u,v ) $ in (\ref{kron}) by  
\begin{eqnarray*}
\CC_{ N, \chi }(X,Y,\t,T) \= F_{\tau}^{\chi}( XT ,  YT )
   F_{\tau}^{\overline\chi}(  {T},  -XYT).
\end{eqnarray*}
\end{theorem}

 Note when $\chi=1$ (so $N=1$) is the trivial character, $\Cg_{f_{1}}(X,Y)=(XY)^{k-2}\Cg_{f_{1}}(-\frac{1}{X}, -\frac{1}{Y})$ and is also $\Rg_f(X,Y)$ and hence $R_{f_{1}}(X,Y)$ is   $R_f(X,Y)$.
At  first glance one might expect that when $\chi$ is trivial we would recover the results of \cite{CPZ} for untwisted modular forms of level $N$, but the requirement that $\chi$ be  primitive means that when it is trivial, $N=1$   (since $\chi(0)=1$)   and we recover   Zagier's original result.
The theorem should be interpreted as a demonstration that together the functions $F_\tau^\chi$ and $F_\tau^{\overline{\chi}}$ encode information specifically about modular forms twisted by characters of conductor equal to their level.

\bigskip
 
  The organization of the paper follows. 
Section 2 is a preliminary section reviewing the theory of  modular forms on $\Gamma_0(N)$ for square free $N$.  We have introduced three different Eisenstein series  twisted    by Dirichlet character. 
In Section 3,  we discuss some properties of period polynomials and explicitly express the period polynomials of various Eisenstein series.
In section 4,  a family of  twisted Kronecker series is introduced and various features of  the twisted series   are described.
 Finally, the proofs of the results are presented in Section 5.

 \medskip

 \noindent  {\bf{Acknowledgment  }}\, \, We would like to thank the referees for numerous helpful comments and suggestions which greatly improved the exposition of this paper.

  \section{\bf{ Preliminary}}
\subsection{\bf{Notations}}
  \begin{itemize}
\item $\mathbb{H}$ : the complex upper half plane
\item $q=e^{2\pi i \tau}, \tau=x+iy \in \mathbb{H}, re(\tau)=x, Im(\tau)=y.$
\item $\G_0(N)=\{\sm a& b\\c& d\esm \in SL_2(\mathbb{Z}) \, : c\equiv 0 \pmod{N} \}, N$ square free
\item $M_{k,N}:=M_k(\Gamma_0(N)),
\,S_{k,N}:=S_k(\Gamma_0(N))\,$  : the space of modular  forms and  cusp forms, respectively,  of even weight $k$ on $\Gamma_0(N) $
\item $\mathcal{B}_{k,N}$ : the basis of Hecke eigenforms in $M_{k,N}$
\item $ \mathcal{B}^{\rm{cusp}}_{k,N} = \mathcal{B}_{k,N} \cap S_{k,N}$
 
\item $W_M\in \sm M\mathbb{Z} & \mathbb{Z} \\ N\mathbb{Z} & M\mathbb{Z} \esm, \det{W_M}=M$ : Atkin-Lehner involution
\item $\j_d=\sm \! d&0\\\!0&1\!\esm$, $W_N=\sm 0 & -1\\ N &   0\esm  $

\item $(f|_k\gamma)(\tau):=\det{\gamma}^{\frac{k}{2}}
(c\tau+d)^{-k}f(\frac{a\tau+b}{c\tau+d}), \gamma =\pm\sm a & b\\ c& d\esm \in PGL_2^+(\mathbb{R})$

\item  $\langle f, g \rangle_N :=\int_{\Gamma_0(N)\backslash \mathbb{H}}f\cdot \overline{g}  \cdot y^{k}\frac{dxdy}{y^2},  \mbox{\, for $f, g\in S_k(\Gamma_0(N))$} $: the Petersson inner product

\item $\chi$ : a primitive   even    character with conductor $N \geq 1  $

\item $\eta=e^{v}, \xi=e^u, u, v \in \mathbb{C}$
\item  $\delta_{*,*}$ : Kronecker Delta-function 
\item  $W( \chi )=\sum_{h\pmod{N}}  
\chi(a) e^{\frac{2\pi i h}{N}}$: the Gauss sum  
\item  $\zeta(s)=\sum_{n\geq 1}\frac{1}{n^s}$ : Riemann Zeta function 
 \item $L({\chi}, s)= 
\sum_{ \tiny{
 n\geq 1,  
\gcd{(N, n)}=1
}} 
\frac{\chi(n)}{n^s}$
\item  $B_k$: the $k$th Bernoulli number     
\item $B_{k,\chi}=N^{k-1}\sum_{h\pmod{N}}\chi(h)B_k(\frac{h}{N})$: the $k$ th  twisted Bernoulli number   
\item    $G_k(\tau)=-\frac{B_k}{2k}+\sum_{n\geq 1}\sigma_{k-1}(n) q^n $: a normalized Eisenstein series on   $\Gamma_0(1) $ 
\item $\mathbb{Q}_f(\chi)$: the field generated by  $a_f(n)$ and $\chi$ over $\mathbb{Q}$ with  $f(\tau)=\sum_{n\geq 1}a_f(n)q^n $
\item  $\mathbb{Q}(\chi)$: the field generated by  $\chi$ over $\mathbb{Q} $  
\end{itemize}

\subsection{Modular forms on $\Gamma_0(N) $  }

In this section we give   the canonical basis of Hecke forms for~$\G_0(N) $ with square free $N.$  See 
 \cite{CPZ}  for the detailed information.   

Our notations for modular forms are standard:  a holomorphic function 
$f:\mathbb{H}^* \rightarrow\C$ 
satisfying  $f|_k\g=f$ for \hbox{all~$\g\in\G_0(N)$}  and having the appropriate growth conditions at the cusps  is a modular form of weight $k$ on $\G_0(N).$  The space $M_{k, N}$ of modular forms has three main decompositions into the space of Eisenstein series  $ M_{k,N}^{\rm Eis} $ and that of cusp forms $S_{k,N},$ into the space of newforms   $M_{k,N}^{\rm new}$    and  that of oldforms  $M_{k,N}^{\rm old}, $   and into Atkin-Lehner involutions: 
$$M_{k,N} = M_{k,N}^{\rm Eis}\bigoplus S_{k,N} ,$$
 $M_{k,N}^{\rm Eis}$ is 
 spanned by  $G_k(d\tau)$ with $ d|N,   N \geq 1  $ if $k\geq 2 $  (see \cite{CPZ} for more details).  
Since all Eisenstein series are oldforms 
we get  
$$M_{k,N}= M_{k,N}^{\rm Eis}\bigoplus_{N_1 | N} \bigoplus_{d| \frac{N}{N_1} } 
 S_{k, N_1 }^{\rm new}\bigl|_k\j_d   
$$
where $S_{k, N_1}^{\rm new}$   is  the space of { newforms} in  $S_{k, N_1}.$
 Let $\BB_{k,N}$ be a basis of $M_{k,N}$ consisting with Hecke forms, which are simultaneous normalized Hecke eigenforms.
The finite set $\BB_{k,N}^{\rm{new}}$ of Hecke forms in $M_{k,N}^{\rm{new}}$ (which are called the newforms) forms 
 a basis of the space   $M_{k,N}^{\rm{new}}.$  
Further decompose     $M_{k,N}$ into   eigenspaces   under the Atkin-Lehner involution :
  let  $\mathfrak{D}(N)$  be   the set of divisors of $N,$ made into a group isomorphic to $(\mathbb{Z}/2\mathbb{Z})^t, t=  | \mathfrak{D}(N)|,$ by multiplication $N_1*N_2=\frac{N_1 N_2}{ (N_1, N_2)^2}.$ 
 Let  $\mathfrak{D}(N)^{\vee}$   be   its dual which is a  group of  characters $\epsilon:\mathfrak{D}(N)\rightarrow\{\pm1 \}.$
Then 
$$M_{k,N} = \bigoplus_{{\epsilon}\in \mathfrak{D}(N)^{\vee}} M_{k, N}^{\epsilon}
= \bigoplus_{{\epsilon}\in \mathfrak{D}(N)^{\vee}}  M_{k,N}^{\rm{old}, \e  } \oplus M_{k,N}^{\rm{new}, \e},$$
where 
$$M_{k,N}^{\e }:= \{   f\in M_{k,N} \, : \,\, f|_kW_M=\epsilon(M) f,\,\, \forall M |N \}.$$  For each decomposition  $N=N_1N_2$ and each   $\epsilon_2\in \mathfrak{D}(N_2)^{\vee},$   there is a   linear map
$$\mathcal{L}_{k, N_2}^{\epsilon_2}: M_{k, N_1}\rightarrow M_{k,N}, \, 
\mathcal{L}_{k, N_2}^{\epsilon_2}(M_{k, N_1}^{\epsilon_1})\subset M_{k, N}^{\epsilon_1\epsilon_2} \,\, \,  (\forall \epsilon_1 \in \mathfrak{D}(N_1)^{\vee})$$
given by 
 \begin{eqnarray}
\label{defLkN}  {\mathcal{L}}_{k,N_2}^{\e_2 }(f) \= f\,{\Bigl|}_{k}\Bigl(\sum_{d|N_2}
\e_2(d) \,\j_d\Bigr)
 \=f\,{\Bigl|}_{k}\Bigl(\sum_{d|N_2}\e_2(d) \,W_d\Bigr).
\end{eqnarray}
Then, using induction on the number of prime factors $t$ of $N$   one can show that  
  
$$ M_{k,N}^{\rm{old}, \e}=\bigoplus_{N=N_1N_2}\mathcal{L}_{k,N_2}^{  \e_2  }(M_{k,N_1}^{{\rm new},   \e_1 }) \, \, (\e\in\mathfrak{D}(N)^{\vee}, 
   \e_j  :=\e|_{\mathfrak{D}(N_j)}, j\in \{1,2\} ).$$

  The space   $M_{k,N}^\e$ has a basis $\BB_{k,N}^\e$ given by
$$ \BB_{k,N}^\e \=\bigcup_{N=N_1N_2}\left\{
 \mathcal{L}_{k,N_2}^{  \e_2   }(f) \, : \, f\in \BB_{k, N_1}^{\rm{new}, \, \e_1  }\right\}= \BB_{k,N}^{\rm{Eis}, \,\e} \cup  \BB_{k,N}^{ \rm{cusp}, \, \e}$$
where  $\BB_{k,N}^{\rm{Eis},\ \e}$ consists of the function $G_{k,N}^{\e}:=\mathcal{L}_{k,N}^{\e}(G_{k }) $ for each  $\e\in \mathfrak{D}(N)^{\vee}$ except $\e=1$ in the case $k=2.$
So  a basis   $\BB_{k,N}$  of $M_{k,N}$  is given by 
 $$\BB_{k,N}:=\bigcup_{\e\in \mathfrak{D}(N)^{\vee}} \BB_{k,N}^\e .$$ 
  
It   is known (see \cite{CPZ})   that   for $N=N_1N_2$,  
if $f\in\BB_{k, N }^\e$ has the form $ {\mathcal{L}}_{k,N_2}^{\e_2}(f_1)$ 
for some $f_1\in\BB_{k, N_1}^{\e_1}$,  then the two scalar products
$\la f,f\ra_N=\la f,f\ra_{\Gamma_0(N)}$  and $\la f_1,f_1\ra_{\Gamma_0(N_1)}$ are related by
\begin{eqnarray*}
\label{PSPRatio}  \la f,f\ra_N \= \la f_1,f_1\ra_{\Gamma_0(N_1)} \,\cdot 
   \prod_{\begin{array}{cc} \tiny{ p|N_2}\\ \mbox{$p$ \tiny{prime}}\end{array}} 2  \,\Bigl(p\+\e_2(p) \,a_f(p)\,p^{1-\frac{k}{2} }\+1\Bigr). \end{eqnarray*}

\medskip

\subsection{ Twisted Eisenstein Series}

In this paper, there are   three   different notions of Eisenstein series twisted by Dirichlet characters which are relevant.
The first type is   
\[(G_{k,N}^{\e})_\chi
= \sum_{n\geq1} \sum_{d \vert n}\chi(d)\,\,\chi(\frac{n}{d})\,\, d^{k-1} q^n \]
with  the associated L-function  (\cite{My}  page $177$)
\begin{eqnarray}\label{dec}
L((G_{k,N}^{\e})_\chi, s)=L(\chi,s)\,\, L(\chi,s-k+1).
\end{eqnarray}
The period polynomial of $(G_{k,N}^{\e})_\chi$ appears in the ${\bf{C}}_{k,N, \chi}(X,Y,\t).$
The second type 
is given in the   double series: 
\[G_{k,\chi}(\tau)
\df  \frac{N^k\Gamma(k)}{2(-2\pi i)^kW(\chi)}
\sum_{(m,n)\neq (0,0)}\frac{\chi(n)}{(mN\tau+n)^k} \]
with  the following
Fourier expansions
\begin{equation}\label{nEisen}
G_{k,\chi}(\tau)=
-\frac{B_{k,\overline{\chi}}}{2k}
+ \sum_{n\geq 1}\sum_{0<d \vert n} \overline{\chi}(d)d^{k-1}e^{2\pi in  {\tau} } .
\end{equation}
 \\
The resulting modular form is in $M_{k}(\Gamma_0(N),  \overline{\chi})$.
Acting by $W_N$, we get  the third type of Eisenstein series  
\[
H_{k,{\chi}}(\tau) := N^{-\frac{k}{2}}W(\chi)(G_{k,\chi}|_kW_N)(\tau)
=\frac{\Gamma(k)}{2(-2\pi i)^k}\sum_{(m,n)\neq (0,0)}\frac{\chi(m)}{(m\tau+n)^k}.
\]
  is in $M_k(\Gamma_0(N), \chi)$   with the  Fourier expansion   (\cite{My}  page $270$)

\begin{equation}\label{hEisen} 
H_{k,\chi}(\tau) 
= \sum_{n\geq 1}\sum_{0<d \vert n}
 { {\chi}}(\frac{n}{d}) d^{k-1}e^{2\pi i n{\tau}}.
\end{equation}

As shown on (\cite{My}  page $177$), 
the L-functions for $G_{k,\chi}$ and $H_{k,\chi}$ have decomposition 
\[
L(G_{k,\chi},s) = \zeta(s)L(\overline{\chi},s-k+1),
\,\,
L(H_{k,\chi},s)=   \zeta(s-k+1) L(\chi,s).
\]

The modular forms $G_{k,\overline{\chi}}+H_{k,\chi}$ and their derivatives appear as Laurent coefficients of the twisted Kronecker theta function $F^{\chi}_\tau(u,v)$ (see section 4).

\medskip

\begin{rmk}
  \begin{enumerate}
\item  The Eisenstein series $G_{k, \chi}(\tau)$ in (\ref{nEisen}) is non-trivial only when $\chi(-1)=(-1)^k.$  

\item   $G_{k,\chi}(\t) =H_{k,\chi}(\t) =G_k(\t),  $ \, when   $N=1$ and $\chi=1.$   
\end{enumerate}
\end{rmk}

\subsubsection{ Twisted Bernoulli Numbers}

The twisted Bernoulli numbers $B_{n,\chi}$ by $\chi$ are given by the generating function
\[
\sum_{a = 0}^{N-1} \frac{\chi(a) e^{au} }{e^{aNu}-1}
= \sum_{n \geq 0} \frac{B_{n+1,\chi}}{(n+1)!}u^n, 
\]
which reduces to the classical Bernoulli numbers $B_n$  when $N=1$.
The generating function for the $B_{n,\chi}$ appears in the constant term of the $q$-expansion of $F_\tau^\chi(u,v)$ in this paper.

  \begin{prop} \label{Ber}  
\begin{enumerate}
\item (page 48 \cite{AIK}) 
  For a primitive even character   $\chi,$  
$B_{2r-1,\chi}=0 , \forall r\in \mathbb{N}$ and
  $B_{0, \chi}= \chi(0)= -2 B_{1, \chi}.$  
\\
 \item   
  There is a  relation between twisted Bernoulli numbers $B_{k,\chi}$ and special values of the $L$-function for $\chi $  when $(-1)^k=\chi(-1)$:
\[
L(\chi,1-k)=-\frac{B_{k,\chi}}{k}=\frac{2W(\chi)N^{k-1}\Gamma(k)}{(2\pi i)^k}L(\overline{\chi},k).
\]
 
\end{enumerate}
\end{prop}

\section{\bf{Period polynomials }}
\subsection{ Period polynomial of cusp forms}  The period polynomial $r_f (X)$   for  $f \in S_{k,N}^{\e}$  
is given by
  
\begin{eqnarray*} \label{Newrfdef} 
 r_f(X) \=  \widetilde f(X,\t) \m \e(N) N^{k/2-1}\,X^{k-2}\,\widetilde f\Bigl(-\frac1{NX},\,-\frac1{N\t}\Bigr)\end{eqnarray*}
for any $\t\in\mathbb{H} $, where $\widetilde f(X,\t)$ is  defined by
\begin{eqnarray*} \label{rtildedef}\widetilde f(X,\t) \= \int_\t^\infty f(z)\,(X-z)^{k-2}\,dz \qquad(\t\in\mathbb{H})\, , 
\end{eqnarray*}
so that 
 
\begin{equation}\label{LPolyCusp}
r_f(X)=\sum_{n=0}^{k-2} (-1)^n \sm k-2\\n \esm r_n(f) X^{k-2-n}
\end{equation}
 with ``periods" $r_n(f)$ defined by
\begin{eqnarray*}\label{rnf}   r_n(f) \df \int_0^\infty f(\tau)\,\tau^n\,d\tau \= \frac{i^{n+1} n!  }{(2\pi  )^{n+1}}\,L(f,n+1) \,\,
 \mbox{\, \, ($0\le n\le k-2)$\,. } 
\end{eqnarray*}
Write
$r_f^{\rm ev}(X)$ and $r_f^{\rm od}(X)$ for the even and odd parts of $r_f(X) 
$ and  was shown \cite{CPZ} 
that 
\begin{eqnarray*} 
 r_f|_{2-k}W_N &=& - r_{f|_kW_N},  \\
  r_{k-2-n}(f) &=& (-1)^{n+1}\e(N) \,  N^{-\frac{k}{2}+1+n} \, r_n(f).\\
\end{eqnarray*}
 Together with~\eqref{defLkN} this implies the relationship,
for $f_1\in M_{k,N_1}^{\rm{new}, \e_1},$
\begin{eqnarray*}\label{relPer}  f=  {\mathcal{L}}_{k,N_2}^{\e_2 }(f_1)\in S_{k,N}^{\e} \;\Rightarrow\; r_f(X)\=\sum_{d|N_2} \e_2(d)d^{1-\frac{k}{2}}
 r_{f_1}(dX)\end{eqnarray*}
between the period polynomial of an oldform and the period polynomial of the newform of lower level from which
it is induced. 
\medskip

 \begin{prop}\label{tfun}
Take $f \in M_{k}(\Gamma_0(N)) $ and let $\chi$ be a primitive character with conductor $N.$
Then 

\begin{enumerate}
\item  $f_{\chi}|_kW_{N^2}
=\chi(-1)\frac{W(\chi)}{W(\overline{\chi})} f_{\overline{\chi}}.$ 

\item $r_{k-2-n}(f_{\chi})=(-1)^{ n+1}\chi(-1)\frac{W(\chi)}{W(\overline{\chi})}N^{ 2n+2-k}r_{ n}(f_{\overline{\chi}})$
\end{enumerate}
 \end{prop}

{\pf} of  Proposition \ref{tfun} : 
\begin{enumerate}
\item    Theorem 4.1 in \cite{AL} implies that Proposition is true when $N$ is prime power.
So it is enough to assume that  $N=Q_1Q_2$ with $Q_i=p_1^{a_1}, \gcd{(p_1,p_2)}=1.  a_i\geq 1, i=1,2.$  Denote $\chi_{M}$ a primitive character  with conductor $M.$ From Proposition 3.4 and Proposition 1.4  in \cite{AL} we get 
\begin{eqnarray*} \label{(0.1)}
   f|_kR_{\chi_{Q_1}}( Q_1)|_kW_{Q_1^2}|_kR_{\chi_{Q_2}}( Q_2)|_kW_{Q_2^2}  
 =\overline{\chi_{Q_2} }(Q_1^2)
W(\overline{\chi_{Q_1}})W(\overline{\chi_{Q_2}}) f_{\chi_{N}}|_kW_{N^2} 
  \end{eqnarray*}
where
$R_{\chi_M}(M):=\sum_{u\pmod{M}}\overline{\chi_M(u)}S_M^u, S_M=\sm M&1\\0&M\esm.$\\

On the other hand  the proof of Theorem 4.1 in \cite{AL} shows that 
\begin{eqnarray*} 
 f|_kR_{\chi_{Q_1}}( Q_1)|W_{Q_1^2}|_kR_{\chi_{Q_2}}( Q_2)|_kW_{Q_2^2}
 =   \,\, \chi_{N}(-1)  \chi_{Q_1}^2(Q_2) W(\chi_{Q_1})    W(\chi_{Q_2}) f_{\overline{\chi_{N}}} 
  \end{eqnarray*}
So we have
\begin{eqnarray*}
 f_{\chi_{N}}|_kW_{N^2}  
=          \chi_{N}(-1) \frac{ \chi_{Q_1}^2(Q_2)}{\overline{\chi_{Q_2} }(Q_1^2)}
\frac{W(\chi_{Q_1})    W(\chi_{Q_2}) } { 
W(\overline{\chi_{Q_1}})W(\overline{\chi_{Q_2}}) }
f_{\overline{\chi_{N}}} 
 = \chi_{N}(-1) \frac{W(\chi_{N}) } { 
W(\overline{\chi_{N}}) }
f_{\overline{\chi_{N}}}
\end{eqnarray*}
 since $W(\chi_{Q_1}\chi_{Q_2})=\chi_{Q_1}(Q_2)\chi_{Q_2}(Q_1)W(\chi_{Q_1})W(\chi_{Q_2}).$
\item  Note, using the result (1), 
\begin{eqnarray*}
&&\int_0^{\infty}(f_{\chi}|_kW_{N^2})(iy)y^{s-1}dy
 = i^k\int_0^{\infty}N^{k-2s}(N^2y)^{s+1-k}f_{\chi}(\frac{i}{N^2y})d\frac{1}{N^2y}\\
&&=i^kN^{k-2s}\int_0^{\infty}y^{k-s-1}f_{\chi}(iy)dy 
 =
 \chi(-1)\frac{W(\chi)}{W(\overline{\chi})}\int_0^{\infty}f_{\overline{\chi}}(iy)y^{s-1}dy.
\end{eqnarray*}
So 
\begin{eqnarray*}
i^kN^{k-2s} \frac{\Gamma(k-s)}{(2\pi)^{k-s}}L(f, {\chi},k-s)=\chi(-1)\frac{W( {\chi })}{W(\overline{\chi})}\frac{\Gamma( s)}{(2\pi)^{ s}}L(f, \overline{\chi}, s).
\end{eqnarray*}
This implies the result. 
\end{enumerate}  

{\qed}
\medskip

\subsection{Rationality}
The following result states the rationality of the twisted period $r_{n }(f_{\chi} )$ in   (\ref{tt}),   which is a part of the extension of Eichler-Shimura-Manin theory:

\begin{thm}\label{Razar}(Razar \cite{Ra}, Shimura \cite{Shimura})
Take  $f(\tau)= \sum_{n\geq 1} a_f(n)e^{2\pi i n \tau}\in S_{k,N}$   which is
  an eigenform for all Hecke operators $T_{\ell},$ prime $\ell, \ell\nmid N.$ \,\,
For any character $\chi$ modulo $r,$  $L(f, \chi,s)$ is entire and there exist non-zero constants $w_f^+ $ and $w_f^-$ depending only on $f$ such that for all integers $j, \,\, 0 \leq j\leq k-2,$
\begin{eqnarray*}
 r_{k-2-j}(f_{\chi} ) 
\in \bigl\{
\begin{array}{cc} \mathbb{Q}_f(\chi) w_f^+ & \mbox{\, if $ (-1)^j=\chi(-1)$}\\
 \mathbb{Q}_f(\chi)w_f^-  & \mbox{\, if $(-1)^j\neq \chi(-1)$} 
\end{array}
\bigr\}.
\end{eqnarray*}
\end{thm}  
 \medskip

\begin{thm} \label{rational} The numbers $w_f^+, w_f^{-} \,  (f \in S_{k, N} $  a normalized Hecke eigenform) \, can be chosen in such a way that 
$w_f^+ w_f^-  = \langle f, \, f \rangle.$
\end{thm}
{\pf} of Theorem \ref{rational} From the result  of Kohnen-Zagier (see Theorem in page $202$ in \cite{KZ}) and that of \cite{Ra} imply the result.
{\qed} 
\medskip

For each  $0\leq m,n \leq k-2$ with $ m-n\equiv 1 \pmod{2},  $ we can choose $w_f^+,  w_f^-$ such that
$$\frac{r_{n}(f_{\chi})\, r_{m }(f)}{w_f^{+}  w_f^{-} }  =\frac{r_{n}(f_{\chi})\, r_{m }(f)}{\langle f, \, f \rangle } \in  \mathbb{Q}_f(\chi). $$

\bigskip
\subsection{  { Period polynomial of twisted  Eisenstein series }} \,

\medskip

 The period ``polynomial"
$r_f(X)$  for $f\in M_{k,1}$  is defined by the same formula~\eqref{Newrfdef} as in the cuspidal case, which
is  independent of the choice of~$\t\in \mathbb{H}$, but with $\widetilde f(X,\t)$ now defined by
\begin{eqnarray*}\widetilde f(X,\t) \,=\, \int_\t^\infty \bigl(f(z)-a_f(0)\bigr)\,(X-z)^{k-2}\,dz\+
      a_f(0)\,\frac{(X-\t)^{k-1}}{k-1}. \end{eqnarray*}

For  Eisenstein series $G_k$, the equation (\ref{LPolyCusp}) can be extended by interpreting ${k-2 \choose s}$ as $\frac{\Gamma(k-1)}{\Gamma(s+1)\Gamma(k-1-s)}$ to get the formula  \cite{Ke}:
\begin{equation*}\label{LPolyGen}
  r_{G_k}(X)=\sum_{n\in \Z}(-1)^n  \Gamma(k-1) \frac{i^{n+1} }{(2\pi)^{n+1}}\lim_{s\rightarrow n}\frac{L(G_k, s+1)}{ \Gamma(k-1-s) } X^{k-2-n}.  
\end{equation*}

The followings are period polynomials   of various Eisenstein series which contribute to Theorem \ref{main}:  
\begin{prop} \label{Kro}
\begin{enumerate}
\item  For $N\geq 1,$  

\begin{eqnarray*}  
r_{(G_{k,N}^{\e})_\chi}(X) &=&   \chi(0) \omega_{G_{k}}^+( X^{k-2}-1) \\
 & + &  \omega_{ G_{k} }^-    {W(\chi)}    N^{1-k}    \sum _{\tiny{\begin{array}{cc}
  r,s\geq  {0} \,  \\
  r+s=k\\
 r,s \, \, \mbox{even}  
\end{array}}} 
\frac{B_{r,\chi}}{r!}
\frac{B_{s,\overline{\chi}}}{ s!} (NX)^{r-1} ,\\
\end{eqnarray*}  
 
where \, $\omega_{G_{k}}^- = -\frac{(k-2)!}{2} $  and    $\omega_{G_{k}}^+ = {(2\pi i)^{1-k}}\zeta(k-1) \omega_{G_{k}}^- .$
\\
\item \cite{CPZ}
Let   $G_{k,N}^{\e}:= {\mathcal{L}}_{k,N}^{\e }(G_k), k>2, \e \in \mathfrak{D}(N)^{\vee} $     and 
$G_{2,N}^-={\mathcal{L}}_{2,N}^{-}(G_2).$ 
So  
\begin{eqnarray*}
&& r_{G_{k,N}^{\e}}(X) = \mathcal{L}_{2-k, N}^{\e}(r_{G_k}(X))=\bigl(r_{G_k} |_{2-k}\sum_{d|N} \e(d) V_d\bigl)(X), \\
 && r^{\rm od}_{G_{k,N}^{\e}}(X) =  \sum_{d|N} \e(d)d^{1-\frac{ k}{2}}  r^{\rm od}_{G_k}(dX),\\
 && r^{\rm ev}_{G_{k,N}^{\e}}(X)=   \omega_{G_k}^+\bigl(\e(N)N^{\frac{k}{2} -1 } X^{k-2}-1   \bigr) 
\prod_{\tiny{\begin{array}{cc} p|N \\ \mbox{$p$\tiny{ prime}}\end{array}}}
\bigl(1 + \e(p) p^{1-\frac{k }{2}}\bigr).
\end{eqnarray*}
\item \cite{CPZ}
\begin{eqnarray*}
&&
\frac{\langle G_{k,N}^{\e}, G_{k,N}^{\e} \rangle}{ \langle G_k, G_k \rangle}
=2^t\prod_{\tiny{\begin{array}{cc} p|N \\ \mbox{$p$  {prime}}\end{array}}}(1+\e(p)  p^{\frac{k}{2}})(1+\e(p) p^{1-\frac{k}{2}}),\\
&& \langle G_k, G_k \rangle =\frac{(k-2)!}{(4\pi)^{k-1}}\frac{B_k}{2k}\zeta(k-1)=-
 \frac{i^{k-1}}{2^{k-1}}\cdot \frac{B_k}{k}\cdot \omega_{G_{k}}^+ .
\end{eqnarray*}
\\
 \end{enumerate}
\end{prop}

{\pf} of {\bf{ Proposition \ref{Kro}}}:\,
 (2) and (3) are discussed in \cite{CPZ}.
(1) can be shown similarly to (1) in \cite{Z1991} using that
\[
r_{(G_{k, N}^{\epsilon})_\chi} (X) = \sum_{n=0}^{k-2}\frac{(-1)^ni^{n+1}\Gamma(k-1)}{(2\pi )^{n+1} }  
 \lim_{s\rightarrow n} \frac{L( (G_{k,N}^\e )_\chi , s+1)}{\Gamma(k-1-s)}   X^{k-2-n}
\]
with the facts (\ref{dec}).
{\qed}
  
\medskip

\section{  Twisted Kronecker series}
 
 Although much of the theory works for odd characters, throughout the rest of the paper we take $\chi$  to be a primitive even character modulo $N>0.$ 

Assume $|q  = e^{2\pi i \tau}   |< \vert \eta  =e^v  \vert, \vert \xi   =e^u   \vert <1,    \tau \in \mathbb{H}, u, v\in \mathbb{C}.  $ 
 Let 
\begin{eqnarray}\label{kron}
F^{\chi}_{\tau}(u,v) =  
\frac{1}{2}\bigl( \sum_{ n\in \mathbb{Z}}  \chi(n) \frac{ \eta^n  }{ \xi   q^{  n}-1}
+  \sum_{ m \in \mathbb{Z}}  \chi(m) \frac{ \xi^{ m} }{ \eta q^{  m}-1}\bigr).
\end{eqnarray}
\\
This twisted Kronecker series shares many properties with the untwisted series (cf Section 3 of \cite{Z1991}).

\noindent   

\begin{prop} \label{twist} 
 Fix a positive integer $N\geq 1.$  The   following    hold: 
\begin{enumerate}
  
\item 
\begin{eqnarray*}
   F^{\chi}_{\tau}(u,v)
&=& \frac{1}{2} \sum_{ h=0 }^{ N  } \chi( h)
\bigl(\frac{\eta^h}{\eta^N-1}+\frac{\xi^h}{\xi^N-1}\bigr)
\\ 
&-&\frac{1}{2}\sum_{m,n\geq 1} (\chi(n)+\chi(m))\bigl(\xi^m\eta^n-\chi(-1)\xi^{-m}\eta^{-n}\bigr ) q^{mn}
\end{eqnarray*}
 
 \item\begin{eqnarray*}
F^{\chi}_{\tau}(u,v)&=&\frac{1}{2 W(\overline{\chi})} \sum_{h=0}^{N-1} \overline{\chi}{(h)}\bigl(F_{\tau}(u+2\pi i \frac{h}{N},v)+F_{\tau}(u, v+2\pi i \frac{h}{N})\bigr)
 \end{eqnarray*}
 
\item   (analytic/meromorphic continuation) $F^{\chi}_{\tau}(u,v) $ can be analytically continued to all $\tau \in \mathbb{H}$ and   meromorphically continued for all   $u,v\in \mathbb{C}$.
   The poles of $F^{\chi}_{\tau}(u,v)$ are at $u = \frac{2\pi i r}{N} + 2\pi i n \tau$ when $(r, N) = 1$
and at $u = 2 \pi i r + 2 \pi i n \tau$ when $(n, N)=1$,
as well as at $v = \frac{2\pi i s}{N} + 2\pi i m \tau$ when $(s, N) = 1$
and at $v = 2 \pi i s + 2 \pi i m \tau$ when $(m,N)=1$. 
 
\item (Laurent expansion) 
\begin{eqnarray*}
&&F^{\chi}_{\tau}(u,v)=  \chi(0) \bigl(
\frac{1}{u}+\frac{1}{v} \bigr)   \\
&& -  \sum_{\tiny{\begin{array}{cc}
r,s\geq 0\\
 r+s \mbox{\, odd}
\end{array}}}\bigl(\frac{1}{2\pi i}\frac{d}{d\tau}\bigr)^{min \{r,s\}}
\bigl(G_{|r-s|+1,  \overline{\chi} }( \tau)
+ H_{|r-s|+1,  \chi}(\tau) \bigl)\frac{u^r}{r!} \frac{v^s}{s!}.
\end{eqnarray*}
\item (Fourier expansion)
\begin{eqnarray*}
F_\tau^\chi(u,v)
&=& \frac{1}{2} \sum_{a=0}^{N}\chi(a)\left(\frac{\xi^a}{\xi^N-1}+\frac{\eta^a}{\eta^N-1} \right)\\
&-&  \sum_{n\geq 1} \sum_{d \vert n}\big(\chi(d)+\chi(\frac{n}{d})\big)\sinh(du+\frac{n}{d}v)q^{n}
\end{eqnarray*}
\item (elliptic property) For  any $ m,n  \in \mathbb{Z},$
\[F^\chi_\tau(u+2\pi i(nN\tau+s),v+2\pi i(mN\tau+r)) = q^{-N^2mn}\xi^{-Nm}\eta^{-Nn}F^\chi_\tau(u,v).\]

\item (modular property) 
$$F^{\chi}_{\frac{a\tau+b}{c\tau+d}}(\frac{u}{c\tau+d}, \frac{v}{c\tau+d})
=\chi(d) (c\tau+d)  e^{ \frac{cuv}{ 2\pi i(c\tau+d) }}F^{\chi}_{\tau}(u,v), \forall \sm a&b\\c&d\esm  \in \Gamma_0(N). $$

\end{enumerate}
\end{prop}
 \medskip

  \begin{rmk}\label{tri} 
  \begin{enumerate}
\item In the case when   $N=1,$  
$ F^{\chi=1}_{\tau}(u,v)$ in  (\ref{kron})   equals
$F_{\t}(u,v).$ 
\item    Note that either $\chi(0)=0 $ if $N>1 $  or $\chi(0)=N=1$. 
\end{enumerate}
 \end{rmk}
 
\section{\bf{Proofs}}
 
\subsection{Proof of Proposition \ref{twist}}


 \begin{enumerate}
\item 
 
Take the first of the sum   of $F_{\t}^\chi$ in (\ref{kron}) and separate it into the cases when $n$ is positive, zero, and negative to get
\begin{eqnarray}\label{SubSum1}
\,\,\, && \frac{1}{2}  \sum_{ n\in \mathbb{Z}}   \chi(n)  \frac{ \eta^n  }{ \xi   q^{  n}-1} \nonumber  \\
 &&= \frac{\chi(0)}{2} \frac{ 1  }{ \xi  -1}   + \frac{1}{2}\sum_{ n>0}  \chi(n) \frac{ \eta^n  }{ \xi   q^{  n}-1}
+\frac{\chi(-1)}{2}\sum_{ n>0}  \chi(n) \frac{ \eta^{-n}\xi^{-1}q^n  }{ 1-\xi^{-1}q^n}.  
\end{eqnarray}
The fist summand in (\ref{SubSum1}) is non-zero only when $\chi$ is trivial and $N=1$, so we have the equality
\[
\frac{\chi(0)}{2} \frac{ 1  }{ \xi  -1}
=\frac{\chi(0)}{2} \frac{ \xi^0  }{ \xi^N  -1}.
\]
The second summand  in (\ref{SubSum1}) can be rewritten via a geometric series to get
\begin{align*}
\frac{1}{2}\sum_{ n>0}  \chi(n) \frac{ \eta^n  }{ \xi   q^{  n}-1}
&=\frac{1}{2}\sum_{ a=1}^{N}  \chi(a)\frac{\eta^a}{ \eta^N-1}
-\frac{1}{2}\sum_{ n>0}\sum_{m > 0}  \chi(n) \xi^m \eta^n  q^{  mn}.
\end{align*}
The third summand in (\ref{SubSum1}) can similarly be rewritten (but without the $m=0$ terms)
\[
\frac{\chi(-1)}{2}\sum_{ n>0}  \chi(n) \frac{ \eta^{-n}\xi^{-1}q^n  }{ 1-\xi^{-1}q^n}
=\frac{\chi(-1)}{2}\sum_{ n>0}\sum_{m>0}  \chi(n) \xi^{-m} \eta^{-n}q^{mn}.
\]

Doing the same for the second sum in (\ref{kron}) and combining with the above yields the result.
 
\item Putting   the following identity
$$\chi(n)=\frac{1}{W(\overline{\chi})}\sum_{a=0}^{N-1}\overline{\chi}(a)e^{2\pi i \frac{a}{N}n}.$$
into the definition of $F^{\chi}_{\tau}(u,v)$ we derive the result.
\item
Breaking out the first $MN$ terms of the sums in $n$ of (1), we get
\begin{eqnarray*}
&& 2 F_\tau^\chi(u,v) 
= \sum_{a=0}^{N }  \frac{\chi(a)\xi^a}{\xi^N-1}
+ \sum_{b=0}^{N}  \frac{\chi(b)\eta^b}{\eta^N-1}\\
&&-\sum_{n=1}^{NM}\left(
\eta^n\sum_{i=1}^N \frac{\chi(i)\xi^i q^{ni}}{1-\xi^N q^{nN}}
+\frac{\chi(n)\eta^n\xi q^n}{1-\xi q^n}
-\eta^{-n}\sum_{i=1}^{N}\frac{\chi(-i)\xi^{-i} q^{ni}}{1-\xi^{-N} q^{nN}}
-\frac{\chi(-n)\eta^{-n}\xi^{-1}q^n}{1-\xi^{-1} q^n}
\right)\label{NMform}\\
&&-\sum_{m\geq 1}  \Bigg(
\xi^m\eta^{NM}\sum_{i=1}^{N}\frac{\chi(i)\eta^i q^{mi}}{1-\eta^N q^{mN}}
+\frac{\chi(m)\xi^m\eta^{NM+1}q^m}{1-\eta q^m}\\
&& -\xi^{-m}\eta^{-NM}\sum_{i=1}^{N}\frac{\chi(-i)\eta^{-i} q^{mi}}{1-\eta^{-N} q^{mN}}
-\frac{\chi(-m)\xi^{-m}\eta^{-(NM+1)}q^m}{1-\eta^{-1} q^m}
\Bigg)q^{mNM},
\end{eqnarray*}
 \\
where the infinite series converges when $\lvert re(u)\rvert < 2\pi (NM+1) Im(\tau)$ so long as $Nv \neq \pm 2\pi i N \tau$.
A corresponding formula holds for $\lvert re(v)\rvert < 2\pi (NM+1)Im(\tau)$ by collecting terms in $m$.
The possible singularities are  those arising from the finite sums of the above equation.
The sum $\sum_{a=0}^{N} \frac{\chi(a)\xi^a}{\xi^N-1}$ has possible poles only 
when $u = \frac{2 \pi i r}{N}, r \in \mathbb{Z},$  and 
there are actual poles when $(r,N)=1$.
Similarly, checking the possible poles from
the terms $\eta^n\sum_{i=0}^{N-1}\frac{\chi(i)\xi^i q^{ni}}{1-\xi^N q^{nN}}, $ 
$\frac{\chi(n)\eta^n\xi q^n}{1-\xi q^n}, \frac{\chi(n)\eta^n\xi q^n}{1-\xi q^n}$ and 
$\eta^{-n}\sum_{i=0}^{N-1}\frac{\chi(-i)\xi^{-i} q^{ni}}{1-\xi^{-N} q^{nN}}$  
the locations of poles can be identified.\\

\item   Collect up common powers of $q$ from $(1)$. 
\item

Since 
$ e^{u+v} -  e^{-(u+v)}
=2\sum_{\substack{r+s \equiv 1 \\ r,s \geq 0}}\frac{u^r}{r!}\frac{v^s}{s!} $
   the expression of $F^\chi_\tau(u,v)$ in  $ (1) $ becomes 
\begin{eqnarray*}
&& F^\chi_\tau(u,v) =  
  B_{0,\chi}(\frac{1}{u}+\frac{1}{v})    \\
&&+ \sum_{\tiny{\substack{r+s \, \mbox{\, odd} \\ r,s \geq 0}}}\left(
\frac{B_{r+1,\chi}}{2(r+1)}\delta_{s,0}
+\frac{B_{s+1,\chi}}{2(s+1)}\delta_{r,0}
- \sum_{n=1}^\infty n^{\tiny{\min\{r,s\}}} \sum_{d|n}
\left(\chi(d)+\chi(\frac{n}{d})\right)d^{|r-s|}q^{n}
\right)\frac{u^r}{r!}\frac{v^s}{s!}\\
&&  =  \chi(0)(\frac{1}{u}+\frac{1}{v})              
 - \sum_{\tiny{\substack{r+s \, \mbox{\, odd} \\ r,s \geq 0}}}
\left(\frac{1}{2\pi i} \frac{d}{d\tau}\right)^{\tiny{\min\{r,s\}}}
\left(G_{|r-s|+1,\overline{\chi}}(\tau)+H_{|r-s|+1,\chi} (\tau)\right)
\frac{u^r}{r!}\frac{v^s}{s!}.
\end{eqnarray*}
 
\item   Using $(2)$ and the elliptic property of  $F_\tau(u,v)$   given in \cite{Z1991} we derive  the   result.
 
\item Using the  modular property  of $F_{\t}(u,v)$ (see \cite{Z1991})   
consider a single summand when $N|c$ and $(N,d)=1$

\begin{eqnarray*}
&& 
\overline{\chi}{(\alpha)}F_{\frac{a\tau+b}  {c\tau+d}}\left(\frac{u}{c\tau+d}+2\pi i \frac{\alpha}{N},\frac{v}{c\tau+d}\right) 
\\
&&=\overline{\chi}(\alpha)(c\tau+d)
e^{\frac{c\left(u+2\pi i(c\tau+d)\frac{\alpha}{N}\right)v}{2\pi i(c\tau+d)}}
F_\tau\left(u+2\pi i(c\tau+d) \frac{\alpha}{N},v \right) \\
&&\mbox{($(N,d)=1$, Modular property $F_\tau$) }
\\
&&=\overline{\chi}(\alpha d)\overline{\chi}(d)^{-1}(c\tau+d)
e^{\frac{cuv}{2\pi i(c\tau+d)}}e^{\frac{\alpha c}{N}v}
F_\tau\left(u+2\pi i\frac{\alpha c}{N}\tau + 2\pi i\frac{\alpha d}{N}, v \right) \\
&&\mbox{($N|c$, elliptic property of $F_\tau$) }
\\
&&= \chi(d)(c\tau+d)e^{\frac{cuv}{2\pi i(c\tau+d)}}\overline{\chi}(\alpha d)F_\tau\left(u+2\pi i\frac{\alpha d}{N}, v \right).
\end{eqnarray*}
Similarly,
\begin{eqnarray*}
&&
\overline{\chi}{(\alpha)}F_{\frac{a\tau+b}{c\tau+d}}
\left(\frac{u}{c\tau+d},\frac{v}{c\tau+d} +2\pi i \frac{\alpha}{N}\right)
\\
&&= \chi(d)(c\tau+d)e^{\frac{cuv}{2\pi i(c\tau+d)}}\overline{\chi}(\alpha d)F_\tau\left(u,v+2\pi i\frac{\alpha d}{N} \right). 
 \end{eqnarray*}
 Summing everything together over $1\leq\alpha\leq N-1$ and dividing by $W(\overline{\chi})$ yields 
\begin{eqnarray*}
 F^\chi_{\frac{a\tau+b}{c\tau+d}}\left(\frac{u}{c\tau+d}, \frac{v}{c\tau+d}\right)  
 =\chi(d)(c\tau+d)e^{\frac{cuv}{2\pi i(c\tau+d)}}
F^\chi_\tau(u,v)
\end{eqnarray*}
\noindent since $\alpha d$ ranges over all of the unit residue classes modulo $N$.

\end{enumerate}
\qed
 
 \medskip

\subsection{\bf{Proof of Theorem \ref{main}}}

To prove the main theorem we need several steps.
  We introduce the notation ${\bf{B}}_{N,\chi}(X,Y, \tau, T)$ to be the product of twisted Kronecker series and want to show that ${\bf{B}}_{N,\chi}(X,Y,\tau, T)={\bf{C}}_{N,\chi}(X,Y,\tau, T)$.
We do this by noting that the coefficient of each $T^{k-2}$ in ${\bf{B}}_{N,\chi}(X,Y,\tau, T)$ is a modular form of weight $k$ on $\Gamma_0(N)$ and so it suffices to verify equality of the cuspidal parts and Eisenstein parts separately   to check the identity in Theorem \ref{main}   coefficient by coefficient.

\noindent {\bf{Step 1:}}  Let
\begin{eqnarray*}
&& {\bf{B}}_{N,\chi}(X, Y, \tau, T)= F_{\tau}^{\chi}(XT,YT) F_{\t}^{\overline{\chi}}( T, -XYT) \\
&=&  \chi(0) \frac{(X+Y)(XY-1)}{X^2Y^2T^2}+
\sum_{ k\geq 2 }{ \bf{B}}_{k,N,\chi}(X,Y,\t)T^{k-2}. 
\end{eqnarray*}

For any $ \gamma=\sm a& b\\ c& d \esm$
 in $\GL_2(\R)$, note 
\begin{eqnarray*}{\bf{B}}_{N,\chi}(X, Y,\tau, T)\vert_2 \gamma = \frac{det{(\gamma)}}{(c\tau+d)^2}
\cdot {\bf{B}}_{N,\chi}(X, Y, \frac{a\t+b}{c\t+d}, \frac{\sqrt{\det{\gamma}}\cdot T}{c\tau+d}).\end{eqnarray*}


\noindent{\bf{Step 2:}}{{\bf{(non-cuspidal part)}}}
 The argument in \cite{CPZ} guarantees 
 that it suffices
to show the equality at the   cusps   only.  So one needs to check if
$\bigl({\bf{C}}_{N,\chi}|W_M\bigr)_{\tau \rightarrow i \infty}=\bigl(
{\bf{B}}_{N,\chi}|W_M\bigr)_{\tau \rightarrow i \infty}, $ for all $M|N.$
Using the  same definition as that of   periods of cusp forms given  in (\ref{pp}),   we have  
 
\begin{eqnarray*}\label{CgkN} 
&&
   R_{ (G_{k, N}^{\epsilon})_\chi} (X,Y)
=   \Rg_{ (G_{k,N}^{\e})_\chi }(X,Y)  +   \Rg_{ (G_{k,N}^{\e})_\chi  }(Y,X),  
\\
 &&
  \Rg_{  (G_{k,N}^{\e})_\chi  }(X,Y) 
= \frac{1}{2}\bigl(
\Cg_{ (G_{k,N}^{\e})_\chi  }(X,Y) +
(XY)^{k-2}\Cg_{ (G_{k,N}^{\e})_\chi  }(-\frac{1}{X}, -\frac{1}{Y})\bigr),   
\\
&&
\Cg_{ (G_{k,N}^{\e})_\chi  }(X,Y) =
 \frac{r_{G_{k,N}^\e}^{\rm ev}( {\frac{Y}{N}} )  \, r_{(G_{k,N}^\e)_{\chi} }^{\rm od}(\frac{X}{N}) + r_{(G_{k,N}^\e)_{\chi} }^{\rm ev}({\frac{Y}{N}}) \, r_{G_{k, N}^\e }^{\rm od}(\frac{X}{N}) 
 }{ N^{1-k} 2 (2i)^{k-3} W(\chi)  \langle G_{k,N}^\e, G_{k,N}^\e \rangle}.  
 \end{eqnarray*}
\\
Further we get
 \begin{eqnarray*}
&&  \Rg_{
 (G_{k,N}^\e)_\chi} (X,Y)   \\
&&
= -\big[  \frac{ 
  r_{G_{k,N}^\e}^{\rm ev}({\frac{Y}{N}}) \,\, r_{(G_{k,N}^\e)_\chi}^{\rm od}( {\frac{X}{N}})
+    ( XY)^{k-2}\, r_{G_{k,N}^\e}^{\rm ev}(\frac{-1}{ NY})\,\, r_{  (G_{k,N}^\e)_\chi }^{\rm od}(\frac{-1}{NX}) }{N^{1-k} (2i)^{k-1 }  W(\chi)  \langle G_{k } , G_{k } \rangle   2^t\prod_{p|N}   (1+\epsilon(p) p^{\frac{k}{2}}
)(1+\epsilon(p) p^{1-\frac{k}{2}})}
 \\
&+ &  \frac{  
  r_{ (G_{k,N}^\e)_\chi}^{\rm ev}(  {\frac{Y}{N}}  )\,\, r_{ G_{k, N}^\e  }^{\rm od}( {\frac{X}{N}})
+    ( XY)^{k-2}r_{ (G_{k,N}^\e)_\chi }^{\rm ev}(  \frac{-1}{ NY}) \,\, r_{ G_{k, N}^\e }^{\rm od}(  \frac{-1}{NX}) }{ N^{1-k} (2i)^{k- 1 }  W(\chi)  \langle G_{k } , G_{k } \rangle   2^t\prod_{p|N}   (1+\epsilon(p) p^{\frac{k}{2}}
)(1+\epsilon(p) p^{1-\frac{k}{2}})} \big]
\end{eqnarray*}
\begin{eqnarray*}
 &=& (1+\chi(0) ) \frac{
( \epsilon(N) N^{\frac{ 2-k}{2}} Y^{k-2}-1)
r_{ (G_{k,N}^\e)_\chi}^{\rm od}(  {\frac{X}{N}} )
 + (- Y^{k-2} +\epsilon(N)  N^{\frac{2-k}{2}})
  X^{k-2} \, r_{ (G_{k,N}^\e)_\chi }^{\rm od}(\frac{-1}{ NX})
}
{   -\frac{B_k}{2k}N^{1-k} 2^{t+1} \prod_{p|N} (1+\epsilon(p)p^{\frac{k}{2}})  W(\chi)  
 }\\
 &=&
 (1+\chi(0) ) \frac{
( \epsilon(N) N^{\frac{2-k}{2}} Y^{k-2}-1)\,
r_{ (G_{k,N}^\e)_\chi}^{\rm od}( {\frac{X}{N}} )
+(  \epsilon(N)  N^{\frac{2-k}{2}}-Y^{k-2})
 X^{k-2}\, r_{ (G_{k,N}^\e)_\chi}^{\rm od}(\frac{-1}{NX}) }
{N^{1-k}  2^{t+1}\prod_{p|N}(1+\epsilon(p)p^{\frac{k}{2}})
 G_{k}(i\infty)  W(\chi) }. \\
&&\mbox{Because}\\
&&\bigl( G_{k,N}^{\e}|_kW_M\bigr)(  i  \infty)=  \e(M)
\sum_{d|N}\e(d)d^{\frac{k}{2}} G_k(  i  \infty)=\e(M)
\prod_{\tiny{\begin{array}{cc} p|N\\ \mbox{$p$ \tiny{prime}} \end{array}}}\bigl(1+\e(p) p^{\frac{k}{2}} \bigl) G_k(  i  \infty) ,
\\
&& \mbox{
so we consider}
 \\
&&   \sum_{\e\in  \mathfrak{D}(N)^{\vee}}\frac{\Rg _{ 
(G_{k,N}^\e)_\chi   } (X,Y)}{ (k-2)!} {\bigl(G_{k,N}^\e|_kW_M\bigr) (i\infty) }  
  \end{eqnarray*}
 (we should have added the condition $\e \neq 1$ if $k=2$ to the summation, but this is not necessary since the symmetry property of  $r_f^{\rm ev}(X)=0$ and so $\Rg_{
 f_{\chi}} (X,Y)=0$  for all $f\in M_{2,N}^\e$ if $\e(N)=1$ (see p1388 \cite{CPZ})) 

\begin{eqnarray*}
&=&  \sum_{\e\in  \mathfrak{D}(N)^{\vee}}
\frac{
 N^{k-1} (1+\chi(0) )  }{  2^{t+1} (k-2)!W(\chi)  } 
\bigl[\e(M)(\e( {N} ) N^{  \frac{ 2-k}{2} } Y^{k-2} -1)
r_{(G_{k,N}^\e)_\chi }^{\rm od}(  {\frac{X}{N}} )
\\
&& +\e(M)(- Y^{k-2}+ \epsilon( {N} ) N^{  \frac{2-k}{2} })
X^{k-2}r_{(G_{k,N}^\e)_\chi }^{\rm od}( \frac{-1}{ NX} ) ]  \\
 && \mbox{ (using the fact that $2^{-t}\sum_{\e }\e(d)\e(d')=\delta_{d,d'}$ for $d,d'\in\mathfrak{D}(N)^{\vee}$) }\\
&&=    (1+\chi(0))  \delta_{M,N} 
\sum_{\tiny{
\begin{array}{cc}
r,s\geq 0 \\
  r,s\, \mbox{\, \, even}  \\
r+s=k\\
\end{array}}
}
  \frac{B_{r, \chi}}{2r!}\frac{B_{s,\overline{\chi}}}{2s!}  
\bigl(
 N^{\frac{2-k}{2}}  
   X^{s-1}
- N^{\frac{2-k}{2}} X^{r-1}Y^{k-2} \bigr)   \\
&&   +  (1+\chi(0)) \delta_{M,1} 
\sum_{\tiny{
\begin{array}{cc}
r,s\geq 0 \\
  r,s\, \mbox{\, \, even}  \\
r+s=k\\
\end{array}}
}
 \frac{B_{r, \chi}}{2r!}\frac{B_{s,\overline{\chi}}}{2s!}  \bigl( X^{r-1} -  X^{s-1}Y^{k-2}  \bigr)  
\\
\end{eqnarray*}

For each $	M|N,$
 
\begin{eqnarray*}\label{comp}
&& C_{ N,\chi}  \vert   W_M (X,Y, \infty , T )
=\chi(0)  \frac{(X+Y)(XY-1)}{X^2Y^2T^2}  \\
 &+&   (1+\chi(0)) \delta_{M,1}     \sum_{k\geq 2}  
\sum_{\tiny{
\begin{array}{cc}
r,s\geq 0 \\
r,s\, \mbox{\, \, even}\\
r+s=k\\
\end{array}}}
\frac{B_{r, \chi}}{2r!}\frac{B_{s,\overline{\chi}}}{2s!}
\bigl( X^{r-1}+Y^{r-1} - X^{s-1}Y^{k-2} - X^{k-2}Y^{s-1}  \bigr) T^{k-2}
\nonumber \\
&+& (1+\chi(0))  \delta_{M,N} 
  \sum_{k\geq 2}   \sum_{\tiny{
\begin{array}{cc}
r,s\geq 0 \\
r,s\, \mbox{\, \, even}\\
r+s=k\\
\end{array}}}
N^{\frac{2-k}{2}}\frac{B_{r, \chi}}{2r!}\frac{B_{s,\overline{\chi}}}{2s!}  
\bigl( X^{s-1} + Y^{s-1} - X^{r-1}Y^{k-2} - X^{k-2}Y^{r-1} \bigr) T^{k-2}.
\end{eqnarray*}

Next let us write
\begin{eqnarray}\label{ep}
  F_{\tau}^{\chi}(u,v)
=\sum_{  k>0, m\geq -1   }{g_{k,m,\chi}}(\tau)(u^{k-1}+v^{k-1})(uv)^m,  
\end{eqnarray}

\begin{eqnarray*}
&&  g_{k,m,\chi}= \left\{\begin{array}{ccc}
\frac{-1 }{m! (m+k-1)!}
 \bigl( \frac{1}{ 2\pi i  }   \frac{d }{d\tau }\bigr)^m  (
G_{k, \overline{\chi} }(\tau)+ H_{k,\chi} (\tau)) & \mbox{\, if $k\geq 2, m\geq 0, $}\\
&& \\
    \chi(0)  & \mbox{\,  if $ k=2, m=-1$}  \nonumber \\
  0 & \mbox{\, otherwise}.
\end{array}\right.
\end{eqnarray*}

So that
\begin{eqnarray*}
&&B(X,Y;\tau,T)
= F_\tau^\chi(XT,YT)
F_\tau^{\overline{\chi}}(T,-XYT)  =    \chi(0)\frac{(X+Y)(XY-1)}{X^2Y^2T^2} +
 \\
&& \sum_{\tiny{
\begin{array}{cc} 
k_1,\,k_2\ge  2  \\ 
 k_1+k_2+2m=k\\
  m_1,\,m_2  \ge -1  \\
m_1+m_2=m 
\end{array}}}
 (-1)^{m_2} g_{k_1,m_1, \chi}\,g_{k_2,m_2, \chi}
  (X^{k_1-1}+Y^{k_1-1})(1-(XY)^{k_2-1})(XY)^m T^{k-2}. 
\end{eqnarray*}

  Since 
 $g_{r,m,\chi}\vert_{k}W_M(i \infty)=0$  for
  $m>0,$    we only need to compute the case  when   $m=0 $  for each $r\geq 2. $  
 A computation shows 
 \begin{eqnarray*}
&&\lim_{\tau \rightarrow i\infty}
G_{r,\chi}\vert_{k}W_\divlevel(\tau)
=
\begin{cases}
- \frac{B_{r,\chi}}{2r } & \text{if } \divlevel = 1 \\
0 & \text{otherwise} 
\end{cases}
\\
&&
\lim_{\tau \rightarrow i\infty}
H_{r,\chi}\vert_{r}W_\divlevel(\tau)
=
\begin{cases}
- \frac{W(\chi)}{\level^{\frac{r}{2}}}
\frac{B_{r,\overline{\chi}}}{2r}
& \text{if } \divlevel = \level \\
0 & \text{otherwise}.
\end{cases}
\end{eqnarray*}

Bringing everything together, we have 
\begin{eqnarray*}
\lim_{\tau \rightarrow i\infty}g_{k_1,0,\chi}(\tau)\vert_{k_1}W_{ M } 
&=&\lim_{\tau \rightarrow i\infty}\frac{-1}{k_1 !}\left(G_{k_1,\chi}(\tau)+H_{k_1,\chi}(\tau) \right)\vert_{k_1 }W_{ M }   \\
&=& 
\delta_{M,1}\frac{B_{k_1,\chi}}{2k_1!}
+\delta_{M,N}\frac{W(\chi)}{\level^{\frac{k_1}{2}}}
\frac{B_{k_1,\overline{\chi}}}{2k_1!}.
\end{eqnarray*}

 Using the above together with the symmetry $g_{r,s,\chi}(\tau)=g_{s,r,\chi}(\tau)$ and the fact that $W(\chi)W(\overline{\chi})=\level$, for  each $k\geq 2$ we get  
  
\begin{align*}
&\lim_{\tau \rightarrow i\infty}B_{k,N,\chi}(X,Y;\tau)\vert_k W_{ M}
\\
=&\delta_{M,1} \sum_{\tiny{ \begin{array}{cc} r+s = k \\ r, s \geq 0\\ r, s \, \mbox{\, even}
\end{array}}}
\frac{B_{s,\chi}}{2s!}\frac{B_{r,\overline{\chi}}}{2r!}
\left(- X^{k-2}Y^{r-1} - X^{r-1}Y^{k-2}+X^{s-1}+Y^{s-1}\right)\\
&+\delta_{M,N} \sum_{\tiny{ \begin{array}{cc} r+s = k \\ r, s \geq 0\\ r, s \, \mbox{\, even}
\end{array}}}
\level^{\frac{2-k}{2}}
\frac{B_{s,\overline{\chi}}}{2s!}
\frac{B_{r,\chi}}{2r!}
\left(- X^{k-2}Y^{r-1} - X^{r-1}Y^{k-2} + X^{s-1} +Y^{s-1} \right)\\
&+2\delta_{N,1} \sum_{\tiny{ \begin{array}{cc} r+s = k \\ r, s \geq 0\\ r, s \, \mbox{\, even}
\end{array}}}
\frac{B_{s}}{2s!}
\frac{B_{r}}{2r!}
\left(- X^{k-2}Y^{r-1} - X^{r-1}Y^{k-2} + X^{s-1} +Y^{s-1} \right).
\end{align*}\\

All together we get  
 \begin{align*}
&B_{N,\chi}|W_M(X,Y; i \infty,T) \\
=&  \chi(0)\frac{(X+Y)(XY-1)}{X^2Y^2T^2}\\
&+\delta_{M,1}\sum_{k\geq 2}
\sum_{\tiny{\begin{array}{cc} r+s = k \\ r, s \geq 0\\ r, s \, \mbox{\, even}
\end{array}}}
\frac{B_{s,\chi}}{2s!}\frac{B_{r,\overline{\chi}}}{2r!}
\left(- X^{k-2}Y^{r-1} - X^{r-1}Y^{k-2}+X^{s-1}+Y^{s-1}\right)T^{k-2}\\
&+\delta_{M,N}\sum_{k \geq 2}
\sum_{\tiny{\begin{array}{cc} r+s = k \\ r, s \geq 0\\ r, s \, \mbox{\, even}
\end{array}}}
\level^{\frac{2-k}{2}}
\frac{B_{s,\overline{\chi}}}{2s!}
\frac{B_{r,\chi}}{2r!}
\left(- X^{k-2}Y^{r-1} - X^{r-1}Y^{k-2} + X^{s-1} +Y^{s-1} \right)T^{k-2}\\
&+2\chi(0)\sum_{k\geq 2}
\sum_{\tiny{\begin{array}{cc} r+s = k \\ r, s \geq 0\\ r, s \, \mbox{\, even}
\end{array}}}
\frac{B_{s}}{2s!}
\frac{B_{r}}{2r!}
\left(- X^{k-2}Y^{r-1} - X^{r-1}Y^{k-2} + X^{s-1} +Y^{s-1} \right)T^{k-2}.
\end{align*}
This agrees  with   
 $C_{N,\chi}|W_M)(X,Y;i\infty,T)$ in (\ref{comp})  so that 
  $$  C_{N,\chi}|W_M (X,Y;i\infty,T)=B_{N,\chi}|W_M(X,Y; i \infty,T)$$  
  for every $M \vert N$.\\  

\noindent {\bf{Step 3:}}{\bf{(cuspidal parts)}} 
To check cuspidal parts  it is enough to check that for any Hecke form $f$ in $\BB_{k,N}^{\rm cusp},$ we get equal Petersson inner products $\langle {\bf{C}}_{k,N,\chi},f \rangle = \langle {\bf{B}}_{k,N,\chi},f \rangle$.
To do that, we will start by decomposing ${\bf{B}}_{k,N,\chi}$ into manageable pieces.  

Using the notation in (\ref{ep}) we get 
 \begin{eqnarray*}
&&  {\bf{B}}_{k,N,\chi}(X,Y,\t)  \\
&=& \sum_{  \tiny{\begin{array}{cc}    k_1,\,k_2  \geq  2, 
  m\ge 0 \\  k_1+k_2+2m=k   
\end{array}}  } 
  (X^{k_1-1}+Y^{k_1-1})(1-(XY)^{k_2-1})(XY)^m \,g_{k_1,k_2,m, \chi}(\t)   
\end{eqnarray*}
with
\begin{eqnarray*}
g_{k_1,k_2,m, \chi} \,& := & \sum_{ \substack{ m_1,\,m_2  \ge -1 \\ m_1+m_2=m \geq 0}  } (-1)^{m_2}g_{k_1,m_1, \chi}\,  g_{k_2,m_2, \overline{\chi}}. 
\end{eqnarray*}

For quasimodular forms $f$ and $g$ of weights $k_1$ and $k_2$, respectively, define the $m$th  modified  Rankin-Cohen brackets by (see \cite{Z1991, CPZ} )
\begin{align}\label{RCB1}
[f,g]_m
:=& [f,g]^{(k_1,k_2)}_m \\
& + \frac{\chi(0)}{2 (2\pi i)^{m+1}}\left(\frac{ \delta_{k_2,2}}{  (m+k_1)}
\frac{d^{m+1} }{d\tau^{m+1}}f
+\frac{(-1)^m  \delta_{k_1,2}}{ (m+k_2)}
\frac{d^{m+1}}{d\tau^{m+1}}g
\right)\nonumber
\end{align}
where the traditional Rankin-Cohen bracket on the space of  modular forms  is defined as
\[
[f,g]^{(k_1,k_2)}_m
:= \frac{1}{(2\pi i)^m}  \sum_{\substack{m_1,m_2 \geq 0 \\ m_1+m_2=m} } (-1)^{ m_2 } 
  {k_1+m-1\choose m_2}{k_2+m-1 \choose m_1} 
 \frac{d^{m_1}}{d\tau^{m_1}}  f
 \frac{d^{m_2}}{d\tau^{m_2}} g.
\]
Then
it  turns out (see \cite{Z1991, CPZ} ) that 
\begin{eqnarray*}
g_{k_1,k_2,m, \chi} \, &=&
\frac{ [G_{k_1, \overline{\chi}}+  H_{k_1, \chi } ,\,G_{k_2, \chi} +  H_{k_2, \overline{\chi}}]_m}{(k_1+m-1)!\,(k_2+m-1)!}.
\end{eqnarray*}

 
\noindent When $\chi$ is trivial, it is known \cite{Z}
that $[  f,\,\, g  ]_{m} $ is in $M_{k,N},   k=k_1+k_2+2m, $ for any $m \geq 0,$
  even when $f$ or $g$ are the quasimodular form $G_2$. It is also 
straightforward to check   that this is  
still true when $\chi$ is an even primitive character modulo $N.$

\noindent {\bf{Step 4:}} {\bf{(Rankin-Selberg )}} 
  In order to compute $\langle {\bf{B}}_{k,N,\chi},f \rangle$ for $f$ in $\BB_{k,N}^{\rm cusp},$    we   need the following proposition.

\begin{prop} \label{RC} For $k_1,k_2> 0$ even and $m\geq 0$ the function 
$g_{k_1,k_2,m, \chi} $ is a modular form of weight $k=k_1+k_2+2m$ on $\Gamma_0(N)$, and its Petersson scalar product with any $f \in \mathcal{B}_{k,N}^{\rm cusp}$ is given by
\begin{eqnarray*}
&& (2i)^{k-1} (k-2)!   W(\chi)  \langle g_{k_1,k_2,m, \chi},f  \rangle \\
&&
=  \sm k-2\\ m\esm \sm k-2\\ m+k_1-1\esm 
 N^{k_2}  \bigl(
r_{k-2-m}(f)\,\,r_{k_2+m-1}(f_{\chi})
+ r_{k-2-m}(f_{\chi})\,\, r_{k_2+m-1}(f )\bigr) \\
&&
  -  \sm k-2\\ m\esm \sm k-2\\ m+k_1-1\esm 
N^{2-k_2}\bigl(
r_{m}(f)\,\, r_{k_1+m-1}(f_{\chi})
+ r_{m}(f_{\chi})\,\, r_{k_1+m-1}(f )\bigr)  
\end{eqnarray*} 
\end{prop}
 
 To prove Proposition\ref{RC} we use the following lemmata.
 
\begin{lemma} \label{w}
For any $f\in M^{  \epsilon }_{k,N},$  
\begin{enumerate}
\item \label{co1} \begin{eqnarray*}
\langle [G_{k_1,\overline{\chi}}, H_{k_2,\overline{\chi}} ]_m \vert W_N,f \vert W_N \rangle
=
N^\frac{k_1-k_2}{2}\epsilon(N)
\left\langle \left[ H_{k_1,\overline{\chi}},G_{k_2,\overline{\chi}}\right]_m,f \right\rangle. 
\end{eqnarray*}
\item \label{co2}
\begin{eqnarray*} 
\langle [ H_{k_1,\chi}, H_{k_2,\overline{\chi}} ]_m \vert W_N,f \vert W_N \rangle
=
N^{1-\frac{k_1+k_2}{2}}\epsilon(N)
\left\langle \left[G_{k_1,\chi},G_{k_2,\overline{\chi}}\right]_m,f \right\rangle.\nonumber 
 \end{eqnarray*}
\end{enumerate}
\end{lemma}
{\pf} of Lemma \ref{w}
\begin{enumerate}
\item \begin{eqnarray*}
&&\langle [G_{k_1,\overline{\chi}}, H_{k_2,\overline{\chi}} ]_m \vert W_N,f \vert W_N \rangle
 =\langle [G_{k_1,\overline{\chi}}\vert W_N, H_{k_2,\overline{\chi}} \vert W_N]_m ,f \vert W_N \rangle  \\
&&=\left\langle \left[\frac{N^\frac{k_1}{2}}{W(\overline{\chi})}  H_{k_1,\overline{\chi}},
\frac{W(\overline{\chi})}{N^{\frac{k_2}{2}}}G_{k_2,\overline{\chi}}\right]_m,\, \epsilon(N) f \right\rangle    \nonumber       \\
 && =
N^\frac{k_1-k_2}{2}\epsilon(N)
\left\langle \left[ H_{k_1,\overline{\chi}},G_{k_2,\overline{\chi}}\right]_m,f \right\rangle. \nonumber 
\end{eqnarray*}
\item
\begin{eqnarray*}
&& \langle [ H_{k_1,\chi}, H_{k_2,\overline{\chi}} ]_m \vert W_N,f \vert W_N \rangle
 =\langle [  H_{k_1,\chi} \vert W_N, H_{k_2,\overline{\chi}}  \vert W_N]_m,f \vert W_N \rangle  \\
&& =
N^{1-\frac{k_1+k_2}{2}}\epsilon(N)
\left\langle \left[G_{k_1,\chi},G_{k_2,\overline{\chi}}\right]_m,f \right\rangle.\nonumber 
 \end{eqnarray*}
\end{enumerate}
{\qed}

\begin{lemma} \label{lemm}
For any $f\in B_{k,N},$  
\begin{enumerate}
\item
\begin{eqnarray*}
\langle [G_{k_1 ,  {\overline\chi}}, G_{k_2,  {\chi}}]_{m}, f \rangle =
\frac{N^{  k_2 }}{(2i)^{k-1}  W(\chi) }
\frac{\Gamma(k-1)   }{  m! \Gamma(k-1-m)}r_{k_2+m-1}(f_{\chi} ) \,\, r_{k-2-m}(f),
\end{eqnarray*}
\item
\begin{eqnarray*}
\langle [H_{k_1,  \chi}, G_{k_2,  { {\chi}}}]_{m}, f \rangle =
\frac{N^{  k_2 }}{(2i)^{k-1}  W(\chi) } 
\frac{\Gamma(k-1)  }{  m! \Gamma(k-1-m)}r_{k-2-m}(f_{ {\chi}})\,\, r_{k_2+m-1}(f ) .
\end{eqnarray*}
\item
 \begin{eqnarray*}
 \langle [ H_{k_1, \chi},  H_{k_2,  \overline{\chi}}]_{m}, f \rangle =
-\frac{N^{2-k_2}}{(2i)^{k-1}W(\chi)}
\frac{\Gamma(k-1)   }{m! \Gamma(k-1-m)} \,\, r_{k_1+m-1}(f_{ \chi}) \,\, r_{m}(f ).
\end{eqnarray*}
\item
\begin{eqnarray*}
 \langle [G_{k_1,  {\overline\chi}},  H_{k_2,  \overline{\chi}}]_{m}, f \rangle =
-\frac{N^{2-k_2}}{(2i)^{k-1}W( \chi )}
\frac{\Gamma(k-1)   }{m! \Gamma(k-1-m)} \,\,r_{m}(f_{\chi} ) \,\, r_{k_1+m-1}(f). 
\end{eqnarray*}
\end{enumerate}
\end{lemma}

{\pf} of Lemma \ref{lemm} : \,\,
First note that Rankin-Selberg method (see Lemma 1 in \cite{KM}) 
tells us the formula, for 
 $f\in S_{k,N},  g\in M_{k_1}(\G_0(N),\overline{\chi }), k=k_1+k_2,$ of
\begin{eqnarray*}
\langle g G_{k_2,  {\chi}}, f \rangle :=\int_{\G_0(N) \backslash \mathbb{H}}
 \overline{f} 
g  G_{k_2,\overline{\chi}} 
y^{k-2}dxdy
 =
\frac{\Gamma(k-1) B_{k_2, {\overline\chi}}}
{(4\pi)^{k-1} 2k_2}
L(g,f; k-1)
\end{eqnarray*}
with   $L(g,f ;s)=\sum_{n\geq 1}\frac{\overline{a_f(n)} a_g(n)}{n^s}\, \, (re(s) \gg 1).$
 Moreover, it can be checked  ( see also  Proposition 6 in \cite{Za}) that
\begin{equation*} \label{RS}
\,\, \langle [g, G_{k_2,\chi}]_m, f \rangle = 
  -  \frac{\Gamma(k-1)\Gamma(k_2+m) B_{k_2, {\overline\chi}}}
{  m!   (4\pi)^{k-1} \Gamma(k_2) 2k_2}
\,\, L(g,f; k-m-1)
\end{equation*}
\\

To compute $L(g,f; k-m-1)$ with $g\in \{ H_{k, \overline{\chi}},  G_{k, {\chi}} \}$  we do the following local computation:
the $L$-series  $L(f,s)=\sum_{n\geq 1}\frac{a_f(n)}{n^s}$ of each Hecke form $f$  has an Euler product 
$L(f,s)=\prod_{\ell}L(f, X)_{\ell}, X=  \ell^{-s},$ where the product is over all primes and where each factor $L(f,X)_{\ell}$ is a rational function of $X.$
 
\medskip

Write any Hecke form $f\in \BB_{k,N}^{\e}$ 
as $f=\mathcal{L}_{k,N_2}^{\e_2}(f_1)$ with
 $f_1\in \BB_{k, N_1}^{\rm{new}, \e_1}$ for some decomposition $N=N_1N_2$ and corresponding decomposition   $\e=\e_1\e_2,$   and then $L(f,X)_{\ell}$ is given  (see \cite{CPZ} p 1384)      as:
 
\begin{eqnarray} \label{od}
  L(\mathcal{L}^{\e_2}_{k,N_2}(f_1), X)_{\ell} 
=L(f_1,X)_{\ell} \cdot
\left\{\begin{array}{cc} 1 \, & \mbox{if $\ell \nmid N_2$}\\
(1+\e(\ell) \ell^{\frac{k}{2} }X)  \, & \mbox{ if $\ell| N_2$}
\end{array}
\right\}
\end{eqnarray}
 with 
\begin{eqnarray*}
L(f_1,X)_{\ell}=\left\{\begin{array}{cc} (1-a_{f_1}(\ell)X+\ell^{k-1}X^2)^{-1} \, & \mbox{if $\ell \nmid N_1$}\\
(1+\e_1(\ell) \ell^{\frac{k}{2}-1}X)^{-1} \, &  \mbox{ if $\ell| N_1$}
\end{array}
\right\}.
\end{eqnarray*}

\begin{eqnarray*}
&& L(G_{k_1,\chi},X)_{\ell}
=\frac{1}{(1-X)(1-\ell^{k_1-1} \overline{\chi} (\ell)X)} 
 =\sum_{n\geq1}\frac{\overline{\chi}(\ell)^{n+1}\ell^{(k_1-1)(n+1)}-1}{\ell^{k_1-1}-1} X^n\\
&&
L(H_{k_1,\overline{\chi}}, X)_{\ell}
=\frac{1}{(1- \overline{\chi} (\ell)X)(1-\ell^{k_1-1} X)} 
 =\sum_{n\geq1}
\frac{\ell^{(k_1-1)(n+1)}-\overline{\chi}(\ell)^{n+1}}{\ell^{k_1-1}-\overline{\chi}(\ell) } X^n\\
\end{eqnarray*}

Now we   prove Lemma in detail: 
\begin{enumerate}
\item To show the identity in  (1) of Lemma,  we treat  the case  when $f$ is a newform and  an oldform separately :
\begin{enumerate}
\item  {\bf{(newforms)}} 
 For  $f\in \BB_{k,N}^{\rm{new },\e}$(this is the special case when $N=N_1, N_2=1$) we have 
\begin{eqnarray*}
L(f, X)_{\ell}=
\left\{
\begin{array}{ccll}
\frac{1}{1-a_{f}(\ell)X+\ell^{k-1}X}=\sum_{j\geq 0} \frac{\alpha^{j+1}-\beta^{j+1}}{\alpha-\beta}X^j
& \mbox{ if $\ell \neq N$}\\ 
\frac{1}{1+\e(\ell)\ell^{\frac{k}{2}-1}X} =\sum_{j\geq0}(-\e(\ell)\ell^{\frac{k}{2}-1}X)^j & \mbox{if $\ell |N$}
\end{array}\right\}
\end{eqnarray*} 
with $\alpha+\beta=a_f(\ell), \alpha \beta=\ell^{k-1}.$\\
\noindent (case i) When $\ell\nmid N,$ we have
\begin{eqnarray*}
&& L(G_{k_1, \chi}, f ; X)_{\ell}=
  \sum_{n\geq   0  }\frac{\overline{\chi}(\ell)^{n+1}\ell^{(k_1-1)(n+1)}-1}{\ell^{k_1-1}-1}  \frac{\alpha^{n+1}-\beta^{n+1}}{\alpha-\beta}X^n\\
&&=\frac{ 
(1-\ell^{k+k_1-2} \overline{\chi} (\ell)X^2)}{
(1-\alpha  \ell^{k_1-1} \overline{\chi} (\ell)X)
(1-\beta  \ell^{k_1-1} \overline{\chi} (\ell)X)
(1-\alpha X) (1-\beta X)}\\
&&=\frac{L(f,X)_{\ell} \,\, L(f,  \overline{\chi} , \ell^{k_1-1}X)_{\ell}}{L( \overline{\chi} , \ell^{k+k_1-2}X^2)_{\ell}} .
\end{eqnarray*}
  
\noindent (case ii)  When $\ell|N$ the computation is similar but more simple and still gets
\[L(G_{k_1, \chi}, f ; X)_{\ell }
=\frac{L(f,X)_{\ell } \,\, L(f,  \overline{\chi} , {\ell}^{k_1-1}X)_{\ell }}{L( \overline{\chi} , \ell^{k+k_1-2}X^2)_{\ell}}.\] 
 \medskip

 Similarly an explicit computation shows that
\begin{eqnarray*}
L( H_{k_1, \overline{\chi}},  f ; X)_{\ell} =\frac{L(f, \overline{\chi} , X)_{\ell}\,\, L(f, \ell^{k_1-1}X)_{\ell}}{L( \overline{\chi} , \ell^{k+k_1-2}X^2)_{\ell}} 
\end{eqnarray*}
using  the fact  $ \frac{B_{k_2, { \overline \chi}}}{2k_2 L(\chi, k_2 )}
= - \frac{N^{  k_2 }\Gamma(k_2)}{(2\pi i)^{k_2}  W(\chi) }.$ 
\medskip 

 Summing over all $\e \in \mathfrak{D}(N)^{\vee}$ together with the above computations  the Petersson scalar product of $f\in \mathcal{B}_{k,N}^{\rm new}$ with the (modified) Rankin-Cohen bracket in (\ref{RCB1}) is given as formulas in Lemma \ref{lemm}.

\item  {\bf{(oldforms)}} 
Take $ f=\mathcal{L}^{  \e_2  }_{k,N_2}(f_1)\in \mathcal{B}_{k,N_2}^{ \rm{ old}, \e_2}$ with $f_1\in \mathcal{B}^{\rm{new}, \e_1}_{k,N_1},  N=N_1N_2, N_2>1.$   Then $L(f,X)_{\ell}$ can be computed (see p 1384 in \cite{CPZ}) as follows :  \,
 
   \begin{eqnarray}\label{lo1}
 && L(\mathcal{L}^{\e_2}_{k,N_2}(f_1),X)_{\ell}=\left\{
\begin{array}{cc}
 {(1+\e_2(\ell) \ell^{\frac{k}{2}}X )} L(f_1, X)_{\ell} \, & \mbox{ if $\ell |N_2$},\\
L(f_1, X)_{\ell}, & \mbox{ if $\ell \nmid N_2$},
\end{array}\right\},  \\
 && L(G_{k_1,\chi},X)_{\ell}=\left\{
\begin{array}{cc}
\frac{1}{  1- X  }  \, & \mbox{ if $ \ell|N$} \nonumber \\
\frac{1}{  (1- X)(1-  \overline{\chi} (\ell)\ell^{k_1-1}X)  }, & \mbox{ if $ \ell \nmid N$}
\end{array}\right\}, \nonumber \\
 && L( H_{k_1,\overline{\chi}}, X)_{\ell}=\left\{
\begin{array}{cc}
\frac{1}{  1-\ell^{k_1-1} X  }  \, & \mbox{ if $ \ell|N$}\\
\frac{1}{  (1-  \overline{\chi} (\ell)X)(1- \ell^{k_1-1}X)  }, & \mbox{ if $ \ell\nmid N$} 
\end{array}\right\}.\nonumber 
\end{eqnarray}
 \\
 So, the convolution   $L$-series of $\mathcal{L}^{\e_2}_{k,N_2}(f_1)$ and $G_{k_1, \chi}$  can be computed as 
\begin{eqnarray*}
&&
L(G_{k_1, \chi},   \mathcal{L}^{\e_2}_{k,N_2}(f_1); s) \\
&&=
\prod_{\ell|N_2}(1+\e_2(\ell)  \ell^{\frac{k}{2}}X)_{\ell}\,\, L(f_1,X)_{\ell}
\prod_{\ell\nmid N_2} \frac{L(f_1,X)_{\ell}\,\, L(f,  \overline{\chi} , \ell^{k_1-1}X)_{\ell}}{L( \overline{\chi} , \ell^{k+k_1-2}X^2)_{\ell}}. \\
\end{eqnarray*}
Therefore,    (\ref{od}) together with   the above  computation shows that 
\begin{eqnarray*}
L(G_{k_1, \chi},  f ; s)  =
  \frac{L(f,s)_{\ell}\,\, L(f,  \overline{\chi} , s- k_1+1  )_{\ell}}{L( \overline{\chi} , 2s- k-k_1+2  ) } \mbox{\, \, since $f_{ \overline{\chi} }=(f_1)_{ \overline{\chi} }.$}\\
\end{eqnarray*}
Summing over all $\e \in\mathfrak{D}(N)^{\vee} $ with the above computation  we get, for
 $  f   \in \mathcal{B}_{k,N}^{ old },$ we get the formulas in Lemma \ref{lemm}-(1).\\
\end{enumerate}

\item
Similar computation as above  shows that
\begin{eqnarray*}
&& L(H_{k_1, \overline{\chi}},   \mathcal{L}^{\e}_{k,N}(f_1); s) =
\prod_{\ell|N}L( f, \ell^{k_1-1}X)_{\ell} \prod_{\ell\nmid N}
\frac{L(f_1, \overline{\chi} , X)_{\ell}\,\, L(f_1, \ell^{k_1-1}X)_{\ell}}{L(\overline{\chi}, \ell^{k+k_1-2}X^2)_{\ell}}
\end{eqnarray*}
to get  
 \begin{eqnarray*}
  L(H_{k_1, \overline{\chi}},  f ; s)
 = \frac{L(  f , s-k_1+1  )  \,\, L(  f, \overline{\chi} , s  ) }{L(\overline{\chi}, 2s-k-k_1+2 ) }.\\
\end{eqnarray*}
 Summing over all $\e \in\mathfrak{D}(N)^{\vee} $ with the above computation  we claim, for
 $  f   \in \mathcal{B}_{k,N}^{ old },$ the identity in Lemma \ref{lemm}-(2).
 
 \item   Using the invariance of the inner product under the slash operator $W_N$ as Lemma \ref{w}  and the formulas (\ref{co1}), (\ref{co2})
 we find that, for $f\in \mathcal{B}_{k,N}^{ \e }.$
  \begin{eqnarray*}\label{i}
 && \langle [ H_{k_1,\chi}, H_{k_2,  \overline{\chi}}]_{m}, f \rangle \\
&&=
\frac{N^{\frac{k_2-k_1}{2}}W(\chi) }{(2i)^{k-1}}
\frac{\Gamma(k-1)   }{m! \Gamma(k-1-m)} \epsilon(N)\,\, r_{k_2+m-1}(f_{ \overline{\chi}}) \,\,r_{k-2-m}(f )\end{eqnarray*}

 With   the functional equation in   (\ref{fun1}) and (\ref{tt})  
\begin{eqnarray*}
&&r_{k-n-2}(f) = (-1)^{n+1}   \epsilon(N)   N^{-\frac{k}{2}+n+1 }\,\, r_n(f) \,
  \\
&&
r_{k-2-n}(f_{\overline{\chi}}) = (-1)^{n+1}\frac{W(\overline{\chi})}{W(\chi)}N^{2n+2-k}\,\, r_{n}(f_{\chi})
\end{eqnarray*}
 and    summing over  all $\e \in\mathfrak{D}(N)^{\vee} $ 
we conclude Lemma\ref{lemm}-(3), for  $f\in \mathcal{B}_{k,N}.$  \\

\item Similar computation as in (3)   
\begin{eqnarray*} 
&& \langle [G_{k_1,  {\overline\chi}},  H_{k_2,  \overline{\chi}}]_{m}, f \rangle  \nonumber \\
&&=\frac{N^{\frac{k_2+k_1}{2}}}{(2i)^{k-1}W(\overline{\chi})}
\frac{\Gamma(k-1)   }{m! \Gamma(k-1-m)} \epsilon(N) \,\, r_{k-2-m}(f_{\overline{\chi}} )\,\, r_{k_2+m-1}(f).
\end{eqnarray*}
 
 and    summing over  all $\e \in\mathfrak{D}(N)^{\vee} $
we conclude Lemma \ref{lemm}-(4), for  
$f\in \mathcal{B}_{k,N}.$  
\end{enumerate}
\qed
\medskip

{\pf} of {\bf{ Proposition \ref{RC}}}:\, Proposition \ref{RC} follows immediately from the above Lemma \ref{lemm}.
\bigskip

\noindent {\bf{Step 5 :}} 
Using  Proposition \ref{RC} the scalar product of any $f\in\BB_{k,N}^{\rm{cusp}} $ with~${\bf{B}}_{k,N,\chi }(X,Y,\tau;T)$ is given by
\begin{eqnarray*}
& & \, (k-2)! \bigl\la   {\bf{B}}_{k, N, \chi}(X,Y,\,\cdot\,), f   \big\ra 
  \\
&&=  \frac{ -  N^{k-1} }{ (2i)^{k-1} W(\chi)  } \sum_{\substack{k_1,\,k_2>0,\;m\ge0\\ k_1+k_2+2m=k} } \binom{k-2}m\,\binom{k-2}{m+k_1-1}  
  \bigl(
\frac{r_{k-2-m}(f)\,\,r_{k_2+m-1}(f_\chi ) }{N^{  k-k_2-1  }}
\\
\, \, && + \frac{r_{k-2-m}(f_{\chi})\,\, r_{k_2+m-1}(f ) }{N^{ k-k_2-1  }}
  - \frac{r_{m}(f)\,\, r_{k_1+m-1}(f_{\chi})}{N^{k+k_2-3}}
- \frac{r_{m}(f_{\chi})\,\, r_{k_1+m-1}(f ) }{N^{k+k_2-3}}  
\bigr) \\
&& (X^{k_1+m-1}Y^m +X^mY^{k_1+m-1} 
- X^{k-m-2}Y^{k_2+m-1} - X^{k_2+m-1}   Y^{k-m-2}).
\end{eqnarray*}
  To show details of this computation, 
for ease of reference, we will name the period terms \\

\begin{tabular}{|c|c|c|c|}
\hline
$\frac{r_{k-2-m}(f)\,\,r_{k_2+m-1}(f_\chi ) }{N^{  k-k_2-1  }}$
& $\frac{r_{k-2-m}(f_{\chi})\,\, r_{k_2+m-1}(f ) }{N^{ k-k_2-1  }}$
& $- \frac{r_{m}(f)\,\, r_{k_1+m-1}(f_{\chi})}{N^{k+k_2-3}}$
& $- \frac{r_{m}(f_{\chi})\,\, r_{k_1+m-1}(f ) }{N^{k+k_2-3}}$ \\
\hline
1 & 2 & 3 & 4 \\
\hline
\end{tabular}
\medskip

and we will name the polynomial terms\\

\begin{tabular}{|c|c|c|c|}
\hline
$X^{k_1+m-1}Y^m$
& $X^mY^{k_1+m-1} $
& $- X^{k-m-2}Y^{k_2+m-1}$
& $- X^{k_2+m-1}   Y^{k-m-2}$ \\
\hline
A & B & C & D \\
\hline
\end{tabular}\\

The overall sum can be computed by  considering the sums arising from choices of pairs of period and polynomial terms.
We will specifically compute the 1A and 1D cases in detail and the rest will follow similar way.
Recall    $ k = k_1+k_2 +2m$, where $k_1$ and $k_2$ are positive even   and $m\geq 0.$

\bigskip
\textbf{1A}:
note that $(k_1+m-1) + (m) = k_1+2m-1 = k - k_2 - 1$, so that
\begin{align*}
&\sum_{\substack{k_1,\,k_2>0,\;m\ge0\\ k_1+k_2+2m=k} }
\binom{k-2}m\,\binom{k-2}{m+k_1-1}
\frac{r_{k-2-m}(f)\,\,r_{k_2+m-1}(f_\chi ) }{N^{  k-k_2-1  }}
X^{k_1+m-1}Y^m\\
&=\sum_{\substack{k_1,\,k_2>0,\;m\ge0\\ k_1+k_2+2m=k} }
\binom{k-2}{k_2+m-1}
r_{k_2+m-1}(f_\chi )
\left(\frac{X}{N}\right)^{k_1+m-1}
\binom{k-2}{k-2-m}\,
r_{k-2-m}(f)
\left(\frac{Y}{N}\right)^m
\end{align*}
so each term in this sum is a term from the product $\pm r_{f_\chi}\left(\frac{X}{N}\right)r_{f}\left(\frac{Y}{N}\right)$.
To determine which terms we have, we need to examine what the values $a = k_2+m-1$ and $b = k-2 - m$ take when $k_1$ and $k_2$ vary as $k_1,k_2 \geq 2$ and $k_1+k_2 \leq k$.
Substituting $m=\frac{k-k_1-k_2}{2}$, we have $a=\frac{k-k_1+k_2}{2}-1$ and $b = \frac{k+k_1+k_2}{2}-2$.
Under the range of $k_1$ and $k_2$, we get the values the values $ k-2\geq b$, $ a+b \geq k-2 $, $b>a$, where $a+b$ odd.
Hence we have
\begin{eqnarray*}
 =\sum_{\substack{ a+b \geq k-2, \;odd\\  k-2  \geq b  >a \geq 0} }
\binom{k-2}{a}
r_{a}(f_\chi )
\left(\frac{X}{N}\right)^{k-2-a}
\binom{k-2}{b}
r_{b}(f)
\left(\frac{Y}{N}\right)^{k-2-b}.
\end{eqnarray*}

 For \textbf{2B}, \textbf{3D},\textbf{4C}, the contribution of the corresponding terms is the same sum as in \textbf{1A}, with the range of summation replaced in 
 \textbf{2B} :
by    $a+b  \geq k-2 , a+b$ odd and $k-2\geq a>b \geq 0.$    The range of summation replaced in  
\textbf{3D} : by  
 $k-2\geq  a+b  $, $a+b$ odd and $a>b\geq 0.$
The range of summation replaced in  
\textbf{4C} : by  
   $ k-2\geq a+b  $, $a+b$ odd and $b>a\geq 0.$

\bigskip
\textbf{1A+2B+3D+4C}:
Bringing all of these terms together, we see that all of the terms in each sum have the same form, just over a distinct set of indices.
The total set of indices are those $a$ and $b$ such that $0 \leq a \leq k-2$ and $0 \leq b \leq k-2$ where $a+b$ is odd, hence we have
\begin{align*}
&\sum_{\substack{k_1,\,k_2>0,\;m\ge0\\ k_1+k_2+2m=k} }
\binom{k-2}m\,
\binom{k-2}{m+k_1-1}
(1A+2B+3D+4C)\\
&=-r_{f}\left(\frac{Y}{N} \right)^{\rm ev}r_{f_\chi}\left(\frac{X}{N}\right)^{\rm od}
-r_{f}\left(\frac{Y}{N} \right)^{\rm od}r_{f_\chi}\left(\frac{X}{N}\right)^{\rm ev}
\end{align*}

\bigskip
\textbf{1B+2A+3C+4D}:
This is the same as above except interchanging the coefficients of $X$ and $Y$,
yielding a total sum of
\begin{align*}
&\sum_{\substack{k_1,\,k_2>0,\;m\ge0\\ k_1+k_2+2m=k} }
\binom{k-2}m\,
\binom{k-2}{m+k_1-1}
(1B+2A+3C+4D)\\
&=-r_{f}\left(\frac{X}{N}\right)^{\rm ev}r_{f_\chi}\left(\frac{Y}{N} \right)^{\rm od}
-r_{f}\left(\frac{X}{N}\right)^{\rm od}r_{f_\chi}\left(\frac{Y}{N} \right)^{\rm ev}.
\end{align*}

\bigskip
\textbf{1D}:
Note that $2(k-2)-(k_2+m-1)-(k-m-2)=k - k_2-1$, so that
\begin{align*}
&\sum_{\substack{k_1,\,k_2>0,\;m\ge0\\ k_1+k_2+2m=k} }
\binom{k-2}{m}
\binom{k-2}{m+k_1-1}
\frac{r_{k-2-m}(f)\,\,r_{k_2+m-1}(f_\chi ) }{N^{  k-k_2-1  }}
\left(- X^{k_2+m-1}
Y^{k-m-2}\right)\\
&=-(XY)^{k-2}\sum_{\substack{k_1,\,k_2>0,\;m\ge0\\ k_1+k_2+2m=k} }
\binom{k-2}{k_2+m-1}
\frac{r_{k_2+m-1}(f_\chi )}{(NX)^{(k-2)-(k_2+m-1)}}
\binom{k-2}{k-2-m}
\frac{r_{k-2-m}(f)}{(NY)^{(k-2)-(k-2-m)}}
\end{align*}
so each term in this sum is a term from the product $\pm (XY)^{k-2}r_{f_\chi}\left(\frac{1}{NX}\right)r_{f}\left(\frac{1}{NY}\right)$.
To determine which terms we have, note that the values $a = k_2+m-1$ and $b = k-2 - m$ take the same values as they did for the (1A) sum, yielding
\begin{align*}
&=-(XY)^{k-2}\sum_{\substack{a+b \geq k-2 ,\;odd\\ k-2 \geq b >a} }
\binom{k-2}{a}
r_{a}(f_\chi )
\left(\frac{1}{NX}\right)^{k-2-a}
\binom{k-2}{b}
r_{b}(f)
\left(\frac{1}{NY}\right)^{k-2-b}.
\end{align*}

\bigskip
\textbf{1D+2C+3A+4B}:
The sums for 2C, 3A, and 4B go similarly to the above, but using the fact that $2(k-2)-(k_1+m-1)-(m)=k+k_2-3$ to evaluate the 3A and 4B versions.
Bringing all of the terms together, we get that
\begin{align*}
&\sum_{\substack{k_1,\,k_2>0,\;m\ge0\\ k_1+k_2+2m=k} }
\binom{k-2}m\,
\binom{k-2}{m+k_1-1}
(1D+2C+3A+4B)\\
&=-(XY)^{k-2}\left(r_{f}\left(\frac{-1}{NY} \right)^{\rm ev}r_{f_\chi}\left(\frac{-1}{NX}\right)^{\rm od}
+r_{f}\left(\frac{-1}{NY} \right)^{\rm od}r_{f_\chi}\left(\frac{-1}{NX}\right)^{\rm ev}\right)
\end{align*}

\bigskip
\textbf{1C+2D+3B+4A}:
This just swaps the roles of $X$ and $Y$ from 1D+2C+3A+4B, so we get
\begin{align*}
&\sum_{\substack{k_1,\,k_2>0,\;m\ge0\\ k_1+k_2+2m=k} }
\binom{k-2}m\,
\binom{k-2}{m+k_1-1}
(1C)+(2D)+(3B)+(4A)\\
&=-(XY)^{k-2}\left(r_{f}\left(\frac{-1}{NX}\right)^{\rm ev}r_{f_\chi}\left(\frac{-1}{NY} \right)^{\rm od}
+r_{f}\left(\frac{-1}{NX}\right)^{\rm od}r_{f_\chi}\left(\frac{-1}{NY} \right)^{\rm ev}\right).
\end{align*}

\bigskip
Putting everything together, we have
\[
(k-2)! \bigl\la   {\bf{B}}_{k, N, \chi}(X,Y,\,\cdot\,), f   \big\ra 
 = R_{f_{\chi}}(X,Y)   \langle f,f \rangle ,
\]
where $R_{f_{\chi}}(X,Y)$ is given in (\ref{pp}).
So we get
\begin{eqnarray*}
B_{k,N,\chi}^{\rm{cusp}}(X,Y, \tau)=\frac{1}{(k-2)!} \sum_{f\in \mathcal{B}_{k,N}^{\rm{cusp}}}R_{f_{\chi}}(X,Y)f(\t).
\end{eqnarray*}

\noindent With the   computation   of  Eisenstein series part we finally  complete  a proof of Theorem \ref{main}.
{\qed}
 
\bigskip
 \section{ \bf{Conclusion}}

The Kronecker series $F_{\tau}(u,v)$   studied by Kronecker and Weil  has a vast  range of applications  in various places such as topology, geometry, mathematical physics, quantum field theory,  combinatorics and number theory.  In modern language , $F_{\tau}(u,v)$ is a meromorphic Jacobi form with matrix index.
 In \cite{C-hilbert}, similar series defined over certain totally real number fields were connected to  Hilbert modular forms.     In contrast,  this paper describes a way to create an infinite family of Kronecker series   by  considering   twists by characters   $\chi \pmod{N}$.
  These new series turn out to have connections   with generating functions of Hecke eigenforms  on $\Gamma_0(N)$,   whose coefficients involve  special values of twisted $L$-functions.   It would be interesting to explore the further arithmetical applications such as a connection with $p$-adic properties, geometry (see \cite{BC}).
\bigskip

\bibliographystyle{amsplain}

 \end{document}